\documentclass[14pt]{article}

\usepackage{amsmath}
\usepackage{amssymb}
\newcommand{\spa}{\hspace*{.1cm}}

\newcommand{\f}{\frac}

\DeclareMathOperator{\Aut}{Aut}

\newcommand{\id}{\mathrm{id}}

\DeclareMathOperator{\stab}{stab}

\newcommand{\C}{\mathbf C}
\newcommand{\R}{\mathbf R}
\newcommand{\Q}{\mathbf Q}
\newcommand{\N}{\mathbf N}
\newcommand{\Z}{\mathbf Z}

\newcommand{\om}{\omega}

\newcommand{\mc}{\mathcal}

\DeclareMathOperator{\ext}{ext}

\begin{document}
\title{Distinguishing extension numbers for $\R^n$ and $S^n$}
\author{Alex Lombardi \\ Harvard University, Cambridge, MA \\ \texttt{alexlombardi01@college.harvard.edu}}
\maketitle

\begin{abstract}
In the setting of a group $\Gamma$ acting faithfully on a set $X$, a $k$-coloring $c: X\rightarrow \{1, 2, ..., k\}$ is called $\Gamma$-distinguishing if the only element of $\Gamma$ that fixes $c$ is the identity element. The distinguishing number $D_\Gamma(X)$ is the minimum value of $k$ such that a $\Gamma$-distinguishing $k$-coloring of $X$ exists. Now, fixing $k= D_\Gamma(X)$, a subset $W\subset X$ with trivial pointwise stabilizer satisfies the precoloring extension property $P(W)$ if every precoloring $c: X-W\rightarrow \{1, ..., k\}$ can be extended to a $\Gamma$-distinguishing $k$-coloring of $X$. The distinguishing extension number $\text{ext}_D(X, \Gamma)$ is then defined to be the minimum $n$ such that for all applicable $W\subset X$, $|W|\geq n$ implies that $P(W)$ holds. In this paper, we compute $\text{ext}_D(X, \Gamma)$ in two particular instances: when $X = S^1$ is the unit circle and $\Gamma = \text{Isom}(S^1) = O(2)$ is its isometry group, and when $X = V(C_n)$ is the set of vertices of the cycle of order $n$ and $\Gamma = \Aut(C_n) = D_n$, the dihedral group of a regular $n$-gon. This resolves two conjectures of Ferrara, Gethner, Hartke, Stolee, and Wenger. In the case of $X=\R^2$, we prove that $\text{ext}_D(\R^2, SE(2))<\infty$, which is consistent with (but does not resolve) another conjecture of Ferrara et al. On the other hand, we also prove that for all $n\geq 3$, $\ext_D(S^{n-1}, O(n)) = \infty$, and for all $n\geq 3$, $\ext_D(\R^n, E(n))=\infty$, disproving two other conjectures from the same authors.
\end{abstract}
\section{Introduction}
Let $\Gamma$ be a group which acts faithfully on a set $X$. As defined by Tymoczko in [7], a $k$-coloring $c: X\rightarrow \{1, ..., k\}$ is distinguishing with respect to $\Gamma$ if the only $\gamma\in \Gamma$ for which $c\circ \gamma = c$ is the identity element (that is, no nontrivial action of some $\gamma\in \Gamma$ fixes the coloring). The distinguishing number of $(X, \Gamma)$, denoted $D_\Gamma(X)$, is defined to be the smallest $k$ such that $X$ has a $\Gamma$-distinguishing $k$-coloring. A special case of this introduced by Anderson and Collins in [1] takes $X = V(G)$ to be the vertices of a graph, and $\Gamma = \Aut(G)$ to be the automorphism group of the graph. Of particular interest are the cases $G = C_n$, the cycle of order $n$. In $[1]$, it is proved that $D(C_n) = 2$ for all $n\geq 6$, while $D(C_3) = D(C_4)=D(C_5)=3$. 
\\\\ In [2], Ferrara, Gethner, Hartke, Stolee, and Wenger introduce a refinement to the distinguishing number problem, in the form of extending precolorings. For the rest of the paper, we fix $k= D_\Gamma(X)$.  Then, given a subset $W\subset X$ and any precoloring $c: X-W\rightarrow \{1, ..., k\}$, we can ask if it is possible to extend $c$ to a $\Gamma$-distinguishing coloring $c^*: X\rightarrow \{1, .., k\}$. For convenience, we introduce the following notation.
\\\\ \textbf{Definiton.} For $W\subset X$ such that the pointwise stabilizer $\stab_\Gamma(W)$ is trivial, we define the \textit{precoloring extension property} $P(W)$ as follows: $P(W)$ holds if and only if every precoloring $c: X-W\rightarrow \{1, 2, ..., k\}$ can be extended to a distinguishing $k$-coloring of $X$. 
\\\\ Based on this notion, in [2], the notion of a \textit{distinguishing extension number} is introduced.
\\\\ \textbf{Definition.} The \textit{distinguishing extension number} $\text{ext}_D(X, \Gamma)$ is equal to the smallest value of $n$ such that for all $W\subset X$, if $|W|\geq n$ and $W$ is not pointwise stabilized by any nontrivial $\gamma\in \Gamma$, then $P(W)$ holds.
\\\\ In this paper, we investigate $\ext_D(\R^n, \text{Isom}(\R^n))$ and $\ext_D(S^n, \text{Isom}(S^n))$, where $S^n$ denotes the unit $n$-sphere; for the rest of the paper, we use $O(n)$ to denote $\text{Isom}(S^{n-1})$ and $E(n)$ to denote $\text{Isom}(\R^n)$. 
\\\\ In the first half of this paper, we compute $\text{ext}_D(X, \Gamma)$ in two particular cases. One case consists of the graph setting mentioned above; in particular, for $X = V(C_n)$ and $\Gamma = \Aut(C_n) = D_n$, where $C_n$ is the cycle of order $n$ and $D_n$ is the dihedral group of a regular $n$-gon. As mentioned earlier, we already know that for $n\geq 6$, $C_n$ has distinguishing number equal to $2$. The other case consists of $X = S^1$ and $\Gamma = O(2)$; it is easy to see that $D_{O(2)}(S^1) = 2$ as well. In $[2]$, it was proved that $\text{ext}_D(\R, E(1)) = 4$, and some partial results on $C_n$ and $S^1$ were given.
\\\\ \textbf{Theorem} [2]. \textit{If $n\geq 6$ is not divisible by $2, 3$, or $5$, then $\ext_D(C_n) = 4$. Furtherore, $\ext_D(S^1, O(2))\leq 16$.}
\\\\ The authors of [2] also conjectured the exact values for the extension numbers for $S^1$ and the remaining $C_n$, which I prove correct here.
\\\\ \textbf{Theorem 1.} \textit{Let $n\geq 6$. If $4$ and $5$ do not divide $n$, then $\ext_D(C_n) = 4$. If $4\mid n$ but $5\nmid n$, then $\ext_D(C_n) = 5$. If $5\mid n$, then $\ext_D(C_n) = 6$. Finally, $\ext_D(S^1, O(2)) = 6$.}
\\\\ The proof of this theorem involves considering general subsets $W\subset S^1$ of cardinality equal to $4$ or $5$ and investigating when $P(W)$ holds. In order to prove Theorem 1, we prove a somewhat stronger characterization of when $P(W)$ holds in this situation.
\\\\ \textbf{Proposition.} 
\\ $i)$ \textit{Let $W\subset S^1$ and $|W| = 5$. Then, $P(W)$ holds unless $W$ is the set of vertices of a regular pentagon.}
\\ $ii)$ \textit{Let $W\subset S^1$, $|W|=4$, and suppose that the four orbits of elements of $W$ under the translation of order $5$ are distinct. Then, $P(W)$ holds unless $W$ is the set of vertices of a square.}
\\\\ The proposition tells us that the only obstructions to extending all precolorings of $S^1- \{\text {four points}\}$ are the obstructions due to symmetries of $C_4$ and $C_5$. 
\\\\ In the second half of the paper, we consider what happens in higher dimensions; we present fairly concrete examples to demonstrate that for all $n\geq 3$, $\ext_D(\R^n, E(n))=\infty$ and $\ext_D(S^n, O(n))=\infty$. In fact, we prove a stronger result.
\\\\ \textbf{Theorem 2.} \textit{Let $n\geq 3$, $X = \R^n$ or $S^n$, and $\Gamma = \text{Isom}(X)$. Then, there exist uncountable sets $W\subset X$ which have trivial pointwise stabilizer inside $\Gamma$ but do not satisfy $P(W, \Gamma)$.}
\\\\ The extension number $\ext_D(\R^3, E(3))$ was previously conjectured to be finite; in fact, Theorem 2 provides the first known instances in which $\ext_D(X, \Gamma)$ is infinite. The question in two dimensions is harder to resolve; in the case of $\R^2$, the authors of $[2]$ conjectured the following.
\\\\ \textbf{Conjecture} [2]. $\ext_D(\R^2, E(2))=7$.
\\\\ We obtain some partial results by considering subgroups of $E(2)$. Let $SE(2): \{\vec v\mapsto A \vec v + \vec b,\spa A\in SO(2),\spa \vec b\in \R^2\}$. We prove the following theorem.
\\\\ \textbf{Theorem 3.} \textit{$\ext_D(\R^2, O(2))=7$, and $\ext_D(\R^2, SE(2))<\infty$.}
\\\\ In the case of $S^2$, we are able to show that $\ext_D(S^2, O(3))=\infty$ using an entirely different argument from the arguments in higher dimensions. We prove the following theorem. 
\\\\ \textbf{Theorem 4.} \textit{$P(W, SO(3))$ does not hold for any finite subset of $S^2$. Furthermore, assuming the axiom of choice, $P(W, SO(3))$ does not hold for any countable subset of $S^2$.}
\\\\ Finally, after proving Theorem 4, we discuss a few unanswered questions regarding extending precolorings on $\R^n$ and $S^n$.
\section{Theorem 1: Extending precolorings on $S^1$}
\subsection{Preliminaries}
Theorem 1 is concerned with computing $\ext_D(S^1, O(2))$ and $\ext_D(C_n, \Aut(C_n))$. We note that the lower bound $\text{ext}_D(S^1)\geq 6$ was already proven in [2] (and the appropriate lower bounds for all of the $C_n$ were also proven). This was done by embedding $C_4$ or $C_5$ (as appropriate) into $S^1$ and $C_n$, and observing that two colors are insufficient to distinguish $C_4$ and $C_5$ (so $P(C_4, S^1)$ and $P(C_5, S^1)$ do not hold). Therefore, we only need to show that the appropriate values also serve as upper bounds to the extension numbers.
\\\\ First of all, we can quickly eliminate all dependence on $C_n$ and work entirely over $S^1$ (see [2] for the complete framework). The vertices of $C_n$ can be embedded into $S^1$ by a map $\phi$ which sends $\{1, 2, ..., n\}$ to the $n$th roots of unity. Under this embedding, for any $W\subset C_n$, we have the following fact.
\\\\ \textbf{Fact.} $P(\phi(W), S^1)$ implies $P(W, C_n)$.
\\\\ This holds because the setwise stabilizer of $\phi(V(C_n))$ inside $O(2)$ is canonically isomorphic to $\Aut(C_n)$. As a result, for the rest this section, we will consider subsets $W\subset S^1$, and $P(W)$ will always be taken to be over $S^1$. Furthemore, we have the following reduction.
\\\\ \textbf{Observation.} For all $\gamma\in O(2)$, $P(\gamma W)$ holds if and only if $P(W)$ holds.
\\\\ This is true because if the coloring $c: \gamma W\rightarrow \{R, B\}$ is preserved by $\gamma^\prime \in O(2)$, then the coloring $c\circ \gamma: W\rightarrow \{R, B\}$ is preserved by $\gamma^{-1}\gamma^\prime \gamma$.
\\\\ This reduction will be used extensively throughout the rest of the paper in the following way: we say that two subsets $W$ and $W^\prime$ of $S^1$ are $O(2)$-equivalent (written $W\cong W^\prime$) if $W^\prime = \gamma W$ for some $\gamma\in O(2)$. This observation tells us that if $W\cong W^\prime$, then we can always replace $W$ with $W^\prime$ without loss of generality, in order to determine if $P(W)$ holds or not.
\\\\ Finally, to establish notation for the rest of the proof, we identify $S^1 \cong \R/2\pi \Z$. We use $\sigma$ to denote a translation, with $\sigma_a$ denoting the map $x\mapsto x+a$, and use $\tau$ to denote a reflection, with $\tau_a$ denoting the map $x\mapsto -x + 2a$. If $c$ is a $2$-coloring, then we will use $c_-$ to denote the opposite coloring to $c$ (i.e., the unique coloring such that $c(x)\neq c_-(x)$ wherever $c$ is defined).
\subsection{An extension of [2] Theorem 7}
In [2], the following theorem was proved.
\\\\ \textbf{Theorem} \textit{[2, Theorem 7]. Suppose $W\subset S^1$ of cardinality $4$ satisfies the following condition, denoted $T(W)$: the intersection $(W+\f i k)\cap W = \emptyset$ for $2\leq k\leq 5$, $1\leq i\leq k-1$. Then, $P(W)$ holds.}
\\\\ To prove this theorem, the authors prove as a lemma that $T(W)$ implies $R(W)$ where $R(W)$ is the following property: there exists $w_0\in W$ such that $\tau_{w_0}(W-\{w_0\})\cap W=\emptyset$. It is then proved that $T(W)$ and $R(W)$ together imply $P(W)$. The goal of Section $3.2$ is to prove that $R(W)$ alone implies $P(W)$. Later, we will replace condition $T(W)$ with successively weaker translation conditions until we have proven Theorem 1.
\\\\ The proof that $R(W)$ implies $P(W)$ is almost exactly the same as the proof of Theorem 7 in [2]; however, we need to substitute the following lemma for Lemma 4 in [2].
\\\\ \textbf{Lemma 2.2.1.} \textit{Suppose that $W\subset S^1$ of cardinality 4 satisfies $R(W)$, and we have a precoloring $c: S^1-W\rightarrow \{R, B\}$. Then, there are at most six extensions of $c$ to $S^1-\{w_0\}$ which are preserved by either $\tau_{w_0}$ or a translation (called ``forbidden" in [2]).}
\\\\ \textbf{Proof.} Since $\tau_{w_0}(W-\{w_0\})\cap W =\emptyset$, there is at most one extension of $c$ to $S^1-\{w_0\}$ which permits $\tau_{w_0}$. Let $\sigma_{\f 1 2}$ be the translation of order $2$. Then, $\sigma_{\f 1 2}(W)\neq W$, because if equality were to hold, property $R(W)$ would not be satisfied. Therefore, there are at most two extensions of $c$ which are preserved by $\sigma_{\f 1 2}$ (there may be two if $W = \{w_0, a, a+\f 1 2, b\}$ where $b\neq \f 1 2+w_0$). Furthermore, if $c^*$ is an extension of $c$ preserved by $\sigma$ of even order, then it is also preserved by $\sigma_{\f 1 2}$ [either clockwise or counterclockwise iteration of $\sigma$ will avoid crossing $w_0$, and shows us that $\sigma_{\f 1 2}$ will in fact always preserve $c^*$]. 
\\\\ Suppose $c_1$ and $c_2$ are extensions of $c$ which permit $\sigma_1$ and $\sigma_2$ of odd or infinite order, and let $w\in W$ be such that $c_1(w)\neq c_2(w)$. We claim that at least one of $\sigma_1$ and $\sigma_2$ has order $3$. On the contrary, suppose that neither $\sigma_1$ nor $\sigma_2$ had order $3$. If $|\sigma_1| = |\sigma_2|<\infty$, then as in the Lemma 4 argument in [2], we let $\mc O_w$ denote the $\sigma_1$-orbit of $w$. In this situation, we may suppose that $\sigma_1 = \sigma_2$ (as $\sigma_1$ will always be some power of $\sigma_2$). Since $c_1(\mc O_w)\cap c_2(\mc O_w) = \emptyset$, we can conclude that $\mc O_2\subset W$. But $|\sigma_1|>4$ (because $|\sigma|$ is odd); since $|W|=4$, we have a contradiction. Therefore, we may assume that $|\sigma_2 | > |\sigma_1|>3$, so $|\sigma_2|\geq 7$. From here, the argument from [2] applies (it is possible to find an element $x_\mc O$ of any $\sigma_1$-orbit $\mc O$ such that $x_\mc O$ and $\sigma_2(x_\mc O)$ are not in $W$), and we obtain a contradiction. Thus, either $\sigma_1$ or $\sigma_2$ has order $3$. 
\\\\ Now, suppose we have $c_1, c_2$, and $c_3$ (of odd or infinite order) which permit $\sigma_1, \sigma_2$, and $\sigma_3$. By the previous paragraph, we obtain that without loss of generality, $|\sigma_1|=|\sigma_2|=3$, which also means that without loss of generality, $\sigma_1=\sigma_2$. Therefore, $|\sigma_1 (W)\cap W| = 3$, and assuming $c_1\neq c_2$, we conclude that $c_1 = c_2$ outside of $\sigma_1(W)\cap W$, while $c_1 = {c_2}_{-}$ on $\sigma_1(W)\cap W$. If we had a fourth extension $c_4$, we would also have $|\sigma_4|=3$, but then we would obtain that $c_4=c_1$ or $c_4=c_2$, a contradiction. Therefore, there are at most three extensions of $c$ which permit a translation of order greater than two. In total, then, we have at most $3+2+1 = 6$ forbidden extensions, which proves Lemma 3.2. $\hfill\blacksquare$
\\\\ \textbf{Theorem 2.2.2.} \textit{If $W\subset S^1$ of cardinality $4$ satisfies condition $R(W)$, then $P(W)$ holds.}
\\\\ \textbf{Proof.} Suppose that $W\subset S^1$, $|W|=4$, and $R(W)$ holds. Given any precoloring $c: S^1-W\rightarrow \{R, B\}$, Lemma 3.2 tells us that there are at least two non-forbidden extensions of $c$ to $S^1-\{w_0\}$. Let $c^*$ be one such extension, which we may further extend to $S^1$ by choosing a color for $w_0$. Assuming for the sake of contradiction that $c$ cannot be extended to distinguish $O(2)$, Lemma 3.2 tells us that the two colorings $c_R$ (obtained from coloring $w_0$ red) and $c_B$ (obtained from coloring $w_0$ blue) are preserved by reflections $\tau_R, \tau_B$ which do not fix $w_0$. The rest of the proof can be taken almost word for word from [2], with the following caveats:
\\\\ 1) In the proof of Lemma 5 in [2], we know that $w_0+\f 1 2\not\in W$ by $R(W)$. 
\\ 2) In the proof of Lemma 12 in [2], the fact that $|\mc O_0|\geq 6$ (which depends on $T(W)$) is irrelevant to the proof of the lemma, and therefore can be omitted.
\\\\ Otherwise, all arguments carry over exactly as written. $\hfill \blacksquare$
\subsection{A weakening of condition $T(W)$}
In order to prove Theorem 1, we will introduce another translational condition $T^\prime(W)$, which is strictly weaker than $T(W)$, and show that $T^\prime(W)$ implies $P(W)$.
\\\\ \textbf{Condition} $T^\prime(W)$: $(W+\f i k)\cap W=\emptyset$ for $k=4, 5$ and $\gcd(i, k) = 1$.
\\\\ \textbf{Theorem 2.3.1.} \textit{If $W\subset S^1$ of cardinality $4$ satisfies $T^\prime(W)$, then it also satisfies $P(W)$.}
\\\\ \textbf{Corollary 2.3.2.} \textit{$\text{ext}_D(C_n)=4$ for all $n\geq 6$ such that $4\nmid n$, $5\nmid n$.} 
\\\\ Most of the necessary work for Theorem 2.3.1 involves checking that $P(W)$ holds in a few specific cases, which occurs in subsequent lemmas. We will first present a short argument that proves Theorem 2.3.1 assuming those lemmas, and prove the lemmas afterwards.
\\\\ \textbf{Proof of Theorem 2.3.1}. Suppose $W\subset S^1$ of size $4$ satisfies $T^\prime(W)$. By Lemma $2.3.3$, there are four possibilities: either $R(W)$ holds (in which case $P(W)$ holds by Theorem $2.2.2$), $W$ is $O(2)$-equivalent to $\{0, \f 1 2, a, a + \f 1 2\}$ ($a\neq \pm \f 1 4$), or $W$ falls into one of two sporadic cases. Lemmas 2.3.4, 2.3.8, and 2.3.10 show that in each of the latter three cases, $P(W)$ holds, so we are done. $\hfill \blacksquare$
\\\\ \textbf{Lemma 2.3.3.} \textit{Suppose $W\subset S^1$ satisifes condition $T^\prime(W)$. Then, either $W$ satisfies $R(W)$, $W \cong \{0, \f 1 2, a, a + \f 1 2\}$ for some $a\neq \pm \f 1 4\in S^1$, $W\cong \{0, \f 1 3, \f 1 2, \f 2 3\}$, or $W\cong \{0, \f 1 6, \f 1 3, \f 1 2\}$.}
\\\\ \textbf{Proof.} Let $W\subset S^1$ be such that $|W|=4$, $T^\prime(W)$ holds, and $R(W)$ does not hold. Without loss of generality (by application of an automorphism of $S^1$), we may assume that $0\in W$. Since $R(W)$ does not hold, we know that $\tau_0(W-\{0\})\cap W\neq \emptyset$, which means that one of the following two statements is true. 
\\ (1) $\f 1 2\in W$
\\ (2) $\exists a\not\in\{ 0, \f 1 2\}$ such that $\{a, -a\}\subset W$
\\\\ Suppose that $\f 1 2\in W$. In this case, we have $W = \{0, \f 1 2, a, b\}$ for some $a$ and $b$. Then, $\tau_a$ does the following:
$$0\mapsto 2a,$$
$$\f 1 2 \mapsto 2a + \f 1 2,$$
$$a\mapsto a,$$
$$b\mapsto 2a-b.$$ 
Since $R(W)$ does not hold, we know that $\tau_a(W-\{a\})\cap W=\emptyset$. However, we know that $\tau_a$ cannot fix $0$ or $\f 1 2$. If $\tau_a$ fixes $b$, then we have $ b = a+\f 1 2$, as claimed. Otherwise, $\tau_a$ must swap two of $\{0, \f 1 2, b\}$; furthermore, possibly shifting $W$ by $\f 1 2$, we may assume that $\tau_a(0)=\f 1 2$ or $\tau_a(0)=b$.
\\\\ If $\tau_a(0) = \f 1 2$, then we have $a = \f 1 4$ or $\f 3 4$, contradicting $T^\prime(W)$. If $\tau_a(0) = b$, then $b=2a$, and we consider $\tau_b$. Again, we know that $|\tau_b(W)\cap W|\geq 2$, and by the same arguments as for $\tau_a$, we may conclude that $\tau_b$ must swap two of $\{0, \f 1 2, a\}$. Furthermore, $\tau_b$ cannot send $0$ to $\f 1 2$. If $\tau_b(0) = a$, then $4a = 2b = a\rightarrow 3a = 0$, which is one of the exceptions covered by the claim. Finally, if $\tau_b(\f 1 2) = a$, then $4a + \f 1 2 = a\rightarrow 3a = \f 1 2$, which is the last exception covered by the claim. Therefore, if $\f 1 2\in W$, $W$ does fall into one of the listed exceptions.
\\\\ On the other hand, suppose that statement $(1)$ is false; by translational symmetry, we may now assume that $(W+\f 1 2)\cap W = \emptyset$. Furthermore, we know that $W = \{0, a, -a, b\}$ for some $a, b\not\in \{0, \f 1 2\}$. By the assumption that $R(W)$ does not hold, we know that $\tau_a(W-\{a\})\cap W\neq \emptyset$; since we also assumed that $a+\f 1 2\not\in W$, there remain three possibilities: $\tau_a(0) = -a$, $\tau_a(0) = b$, or $\tau_a(-a) = b$.
\\\\ If $2a=\tau_a(0)=-a$, then $3a = 0$ and by symmetry, we may assume that $a = \f 1 3$. Then, $\tau_b$ cannot send one cube root of unity to another unless $b$ is some $6$th root of unity, as included in the list of exceptions. 
\\\\ If $2a = \tau_a(0)=b$, then we consider $\tau_{-a}(W) = \{-2a, -3a, -a, -4a\}$. Since $|\tau_{-a}(W)\cap W|\geq 2$, either $3a, 4a, 5a$, or $6a$ is equal to $0$. But $3a=0\rightarrow b=2a = -a$, a contradiction, while the second two subcases are impossible by $T^\prime(W)$ and the fact that $\f 1 2\not\in W$. The final subcase is one of the listed exceptions.
\\\\ If $3a = \tau_a(-a) = b$, then we consider $\tau_{-a}$; if it does not fall into either of the first two categories, then we obtain the opposite result: $-3a = b$ as well. But then $2b = 0$, a contradiction of Case 2 ($(W+\f 1 2)\cap W = \emptyset$).  Thus, the listed exceptions are in fact the only exceptions, as desired. $\hfill \blacksquare$
\\\\ \textbf{Lemma 2.3.4.} \textit{If $W = \{0, a, \f 1 2, a + \f 1 2\}$ for some $a\neq \pm \f 1 4$, then $P(W)$ holds.}
\\\\ \textbf{Proof.} Since $a\neq \pm \f 1 4$, the collection $\{0, a, \f 1 2, a+\f 1 2, -a, -a+\f 1 2\}$ consists of six distinct points. Let $c$ be any precoloring of $S^1-W$ such that $c(-a)=c(-a+\f 1 2) = R$, and let $d$ be any precoloring such that $d(-a)=R$, $d(-a+\f 1 2) = B$ (by negating colorings, proving Lemma 2.3.4 in these two cases suffices to prove the $(B, R)$ and $(B, B)$ cases as well).
\\\\ Extend $c$ (respectively, $d$) to $c_1$ and $c_2$ (respectively, $d_1$ and $d_2$) in the following way: define $c_2(0)=d_1(0)=R$, $c_1(0)=d_2(0)=B$, $c_i(a)=d_i(a)=B$ for $i=0, 1$, $c_i(\f 1 2)=d_i(\f 1 2) = R$, $c_1(a+\f 1 2)=d_1(a+\f 1 2)=B$, $c_2(a+\f 1 2)=d_2(a+\f 1 2) = R$. 
\\\\ \textbf{Note 2.3.5.} For $k$ equal to $c$ or $d$, $k_1$ and $k_2$ differ only on $W^\prime := \{0, a+\f 1 2\}$.
\\\\ \textbf{Note 2.3.6.} $\tau_0, \tau_{\f 1 4}$, $\tau_{\f {a+\f 1 2}2}$, $\sigma_{\f 1 2}$, $\sigma_a$, and $\sigma_{a+\f 1 2}$ do not preserve any of $c_1, c_2, d_1, d_2$. Furthermore, $\tau_{\f a 2}$ does not preserve $c_1, c_2$, or $d_1$, and $\tau_{-\f a 2}$ does not preserve $c_1, d_1$, or $d_2$. In particular, if $\gamma\in O(2)$ preserves $c_1, c_2, d_1$, or $d_2$, then $\gamma(0)\not\in W^\prime$.
\\\\ Finally, for $k=c$ or $k=d$, let the ``intermediate coloring'' $k_3$ be such that $k_3=k_1=k_2$ on $S^1-\{0, a+\f 1 2\}$, $k_3(0)=k_-(2a)$ (which is well-defined because $2a\not\in \{0, a+\f 1 2\}$), and $d_3(a+\f 1 2)= d_3(0)$ while $c_3(a+\f 1 2)= {c_3}_-(0)$.
\\\\ We prove Lemma 2.3.4 by showing that one of $\{c_1, c_2, c_3\}$ is distinguishing, and one of $\{d_1, d_2, d_3\}$ is distinguishing. The arguments for $c$ and $d$ are extremely similar; we will have $k$ denote either $c$ or $d$ and specify at which points the two arguments differ.
\\\\ Assume that neither $k_1$ nor $k_2$ is distinguishing. Then, $k_1$ and $k_2$ must each be invariant under some nontrivial reflection or translation. \\\\ Suppose $k_1$ and $k_2$ are invariant under translations $\sigma_1$ and $\sigma_2$. Then, we will use the fact that $\sigma_1\sigma_2(0)=\sigma_2\sigma_1(0)$ to derive a contradiction. By definition of $\sigma_1$, $k_1(\sigma_1(0))=k_1(0)$, which implies that $k_2(\sigma_1(0))=k_1(0)$ unless $\sigma_1(0)\in W^\prime$; this is not the case by Note 2.3.6.
\\\\ Therefore, $k_2(\sigma_1(0)) = k_1(0)$. This then implies that $k_2(\sigma_2\sigma_1(0))=k_1(0)$, which implies that $k_1(\sigma_2\sigma_1(0))=k_1(0)$ except in the following two situations.
\\\\ $i)$ $\sigma_2\sigma_1(0)=0$: This would mean that $\sigma_1(0)=\sigma_2^{-1}(0)$, which cannot happen because we know that $k_2(\sigma_1(0))=k_1(0)$, but $k_2(\sigma_2^{-1}(0))=k_2(0)\neq k_1(0)$.
\\\\ $ii)$ $\sigma_2\sigma_1(0)=a+\f 1 2$: since $\sigma_1^2 \neq \sigma$ and $\sigma_1^2\neq \id$ (Note 2.3.6), we may suppose that $\sigma_2\sigma_1(0)\neq a+\f 1 2$ by replacing $\sigma_1$ with $\sigma_1^2$ if necessary (it cannot be the case that $\sigma_2\sigma_1(0)=a+\f 1 2 = \sigma_2\sigma_1^2(0)$). 
\\\\ Since we can arrange that $\sigma_2\sigma_1(0)\neq a+\f 1 2$, we can conclude that $k_1(\sigma_2\sigma_1(0))=k_2(\sigma_2\sigma_1(0))=k_1(0)$. Because this is entirely symmetric in $\sigma_1, \sigma_2$, the same argument (except applying $\sigma_2$ first) proves that $k_1(\sigma_2\sigma_1(0))=k_2(0)$, which is a contradiction.
\\\\ Now suppose that $k_1$ is invariant under a translation $\sigma$ and $k_2$ is invariant under a reflection $\tau$. Again, by replacing $\sigma$ with $\sigma^2$ if necessary, we can arrange that $\tau\sigma(0)\neq a+\f 1 2$. We have the relation $\tau\sigma(0) = \sigma^{-1}\tau(0)$, which we will use to derive a contradiction. In fact, the previous argument fully carries over to allow us to conclude that $k_1(\tau\sigma(0))=k_2(\tau\sigma(0)) = k_1(0)$, while $k_1(\sigma^{-1}\tau(0))=k_2(\sigma^{-1}\tau(0))=k_2(0)$. This argument also applies when $k_1$ is invariant under a reflection and $k_2$ is invariant under a translation (this is essentially a Red-Blue color swap). 
\\\\ Therefore, if neither $k_1$ nor $k_2$ is distinguishing, then they must be invariant under reflections $\tau_1$ and $\tau_2$, respectively. For the remainder of the proof, let the translation $\sigma = \tau_2\tau_1$, and let $k$ be defined on $S^1-W^\prime$ (i.e., extend it to $a$ and $\f 1 2$ because $k_1$ and $k_2$ match there). This step is inspired by the argument in $[2]$ but takes it further -- in $[2]$, two reflections are composed in this way in a situation where their corresponding colorings differ at only one point. As a result, the orbits here are more complicated.
\\\\ \textbf{Observation 2.3.7.} $\sigma$ preserves $k$ on $S^1-\{0, \sigma^{-1}(0), \tau_1(0), a+\f 1 2, \sigma^{-1}(a+\f 1 2), \tau_1(a+\f 1 2)\}$. Additionally, $k$ takes specific values as dictated by the following chart.
\begin{center}
  \begin{tabular}{ | l | r | }
    \hline
    $\sigma^{-1}(0)$ & $k_2(0)$ \\ \hline
    $\sigma(0)$ & $k_1(0)$ \\ \hline
    $\tau_1(0)$ & $k_1(0)$ \\ \hline
    $\tau_2(0)$ & $k_2(0)$ \\ \hline
    $\tau_1(a+\f 1 2)$ & Blue \\ \hline
    $\tau_2(a+\f 1 2)$ & Red \\ \hline
    
    \hline
  \end{tabular}
\end{center}

\noindent \textbf{Justifications.} 1. For $\theta\in S^1$, if $\theta\not\in W^\prime$ and $\sigma(\theta)\not\in W^\prime$, then $k(\theta)=k(\sigma(\theta))$ (which is $k(\tau_2\tau_1(\theta))$) unless $\tau_1(\theta)\in W^\prime$. The colors of the $\tau_i(0), \tau_i(a+\f 1 2)$ are dictated by the reflections' color-preserving properties.
\\\\ 2 (to show $k(\sigma(0))=k_1(0)$). We know that $k(\sigma(0))=k_1(0)$ unless $\tau_1(0)\in W^\prime$ or $\sigma(0)\in W^\prime$. But $\tau_1(0)\not\in W^\prime$ by Note 2.3.6, and $\sigma(0)\neq 0$ because if the opposite were true, then $\tau_1(0) = \tau_2(0)$, despite the fact that $\tau_1(0)$ and $\tau_2(0)$ must have opposite colors. Finally, if $\sigma(0)=a+\f 1 2$, then $\sigma(\f 1 2) = a$. We know that $k(a)=B$ and $k(\f 1 2)=R$ for both $k=c$ and $k=d$. Therefore, there is a violation of $\sigma$ color-preservation in either case ($\sigma$ taking a red point to a blue point). We already have enough information to know that this happens at most at $\tau_1(0)\mapsto \tau_2(0)$ [for $k=c$, $\text{Red}\mapsto \text{Blue}$ in fact never happens], so we would need that $\tau_1(0)=\f 1 2$, which is impossible by Note $2.3.6$.
\\\\ 3 (to show $k(\sigma^{-1}(0))=k_2(0)$). The argument is similar here. As before, $\tau_2(0)\not\in W^\prime$ by Note 2.3.6, and we have already shown that $\sigma^{-1}(0)\neq 0$. If $\sigma^{-1}(0)=a+\f 1 2$, then then we use the fact that then $\sigma(0)=-a+\f 1 2$. We proved already that $k(\sigma(0))=k_1(0)$, but in both the $k=c$ and $k=d$ cases, $k(-a+\f 1 2)\neq k_1(0)$, a contradiction. Therefore, $k(\sigma^{-1}(0))=k_2(0)$. 
\\\\ To complete the proof of Lemma 2.3.4, we now assume for the sake of contradiction that $k_3$ is also not distinguishing.
\\\\ Again, we can see easily that $k_3$ cannot permit a nontrivial translation $\sigma_3$, using the relation $\tau\sigma_3=\sigma^{-1}_3\tau$ for $\tau \in \{\tau_1, \tau_2\}$. Pick $i\in \{1, 2\}$ such that $k_i(0)\neq k_3(0)$ (there is exactly one such $i$). On one hand, $\sigma_3(0)\not \in W^\prime$ (this can be easily checked), so $k(\sigma_3(0))=k_3(\sigma_3(0))=k_3(0)$ and therefore $k_i(\tau_i\sigma_3(0))=k_3(0)$. On the other hand, we already know that $\tau_i(0)\not\in W^\prime$ (by Note 2.3.6) and so $k_3(\tau_i(0))=k_i(\tau_i(0))=k_i(0)$, and hence $k_3(\sigma_3^{-1}\tau_i(0))=k_i(0)$. Since $k_3$ and $k_i$ differ at only $0$, this means that $\sigma_3\tau_i(0))=0$, i.e., $\tau_i(0)=\sigma_3^{-1}(0)$; this contradicts the fact that $\tau_i(0)$ and $\sigma_3^{-1}(0)$ have opposite colors. 
\\\\ Therefore, we conclude that $k_3$ permits another reflection, $\tau_3$, and because $k_3$ permits neither $\tau_0$ nor $\tau_{a+\f 1 2}$ (we chose $k_3$ specifically so this is the case), we have that $\tau_3\neq \tau_1, \tau_3\neq \tau_2$. Thus, we have two nontrivial translations $\sigma_{31}:= \tau_3\tau_1$ and $\sigma_{23}:=\sigma_2\sigma_3$, satisfying the relation $\sigma_{23}\sigma_{31}=\sigma:= \sigma_{21}$.
\\\\ We will derive a contradiction using the fact that $\sigma_{21}$ is also equal to $\sigma_{31}\sigma_{23}$ (i.e., the translations commute). Let $i\in\{1, 2\}$ be such that $k_i(0)\neq k_3(0)$, and let $j \in \{1, 2\}$ be such that $j\neq i$. Observation $2.3.7$ tells us that $k(\sigma_{ji}(0))=k_i(0)$. On the other hand, $\sigma_{ji}(0)=\sigma_{3i}\sigma_{j3}(0)$, and $k_j(0)=k_3(0)$. Since $k_j$ and $k_3$ differ only at $a+\f 1 2$, this means that $k_3(\sigma_{j3}(0))=k_j(0)$ unless:
\\\\ $1)$ $\tau_3(0)=a+\f 1 2$, i.e., $\tau_3 = \tau_{\f{a+\f 1 2}2}$. This does not hold, as we can easily note that $\tau_{\f{a+\f 1 2}2}(a)=\f 1 2$ implies that this particular reflection does not preserve $k_3$.
\\\\ $2)$ $\tau_j\tau_3(0)=a+\f 1 2$. For $k=d$, this does not hold, because $\tau_3(0)$ has the same color as $0$ under both $k_j$ and $k_3$, while $a+\f 1 2$ has the opposite color under $k_j$. For $k=c$, this does not hold, because then we would have $\tau_j\tau_3(\f 1 2) = a $. Since $k(\f 1 2) = R$ and $k(a) = B$, this means that $\tau_3(\f 1 2) = \tau_j(a) \in W^\prime$. In particular, $c_j(\tau_j(a))=B$ implies that $j=1$, while $c_3(\tau_j(a))=R$ implies that $\tau_3(\f 1 2)=\tau_1(a)=a+\f 1 2$. Then, we get a contradiction from the fact that $\tau_1\tau_3(-a+\f 1 2)=0$, while $\tau_3(-a+\f 1 2)=2a+\f 1 2\not\in W^\prime$ [$\tau_1\tau_3$ sends the red $-a+\f 1 2$ to the $c_1$-blue $0$]. 
\\\\ Therefore, $k_3(\sigma_{j3}(0))=k_j(\sigma_{j3}(0))=k_j(0)$. Since $\sigma_{j3}(0)\neq 0$, this means that $k_i(\sigma_{j3}(0))=k_j(0)$ as well. But this implies that $k_i(\tau_i\sigma_{j3}(0))=k_j(0)$, and then we conclude that $k_i(\sigma_{3i}\sigma_{j3}(0))=k_j(0)$ unless $\sigma_{3i}\sigma_{j3}(0)=0$ (but $\sigma_{3i}\sigma_{j3}=\sigma_{ji}$, which we know is nontrivial) or $\tau_i\sigma_{j3}(0)=0$ (which contradicts the fact that $k_i(\tau_i\sigma_{j3}(0))=k_j(0)\neq k_i(0)$). Thus, $k_i(\sigma_{ji}(0))=k_i(\sigma_{3i}\sigma_{j3}(0))=k_j(0)$, contradicting Observation 2.3.7 (which says that $k(\sigma_{ji}(0))=k_i(0)$). Hence, one of $k_1, k_2$, and $k_3$ is distinguishing. This proves Lemma 2.3.4. $\hfill \blacksquare$
\\\\ \textbf{Lemma 2.3.8.} \textit{If $W=\{0, \f 1 3, \f 1 2, \f 2 3\}$, then $P(W)$ holds.}
\\\\ \textbf{Proof.} The principles behind this proof are the same as those behind Lemma 2.3.4, and many of the same reductions are made.  Let $c$ be any precoloring of $S^1-W$ such that $c(\f 5 6)=c(\f 1 6)=R$, and let $d$ be any precoloring of $S^1-W$ such that $d(\f 5 6)=B$ and $d(\f 1 6)=R$. Let $W^\prime = \{0, \f 2 3\}$, and extend $c$ and $d$ to $c_1, d_1, c_2, d_2, c_3$, and $d_3$ in the following way:
\begin{center}
  \begin{tabular}{ | l || c | c | c || c | c | c |}
    \hline
    & $c_1$ & $c_3$ & $c_2$ & $d_1$ & $d_3$ & $d_2$ \\ \hline
    $0$ & B & R & R & B & B & R \\ \hline
    $\f 1 3$ & B & B & B & R & R & R \\ \hline
    $\f 1 2$ & B & B & B & B & B & B \\ \hline
    $\f 2 3$ & R & R & B & R & B & B \\ \hline
  \end{tabular}
\end{center}
In this proof, $c_1, c_2$, and $c_3$ have the same purpose as they did in the proof of Lemma 2.3.4; that is, we first assume that none of $c_1, c_2$, or $c_3$ is distinguishing, and show that $c_1$ and $c_2$ both permit translations $\tau_1, \tau_2$. Then, we will show that the ``intermediate'' coloring $c_3$ also permits a translation $\tau_3$, and derive a contradiction from the (commutative) relation
$$\tau_2\tau_1 = (\tau_2\tau_3)(\tau_3\tau_1) = (\tau_3\tau_1)(\tau_2\tau_3).$$
First, note the following things about the six colorings defined above.
\\\\ \textbf{Note 2.3.9.} $c_1$ and $d_1$ are $\Aut(C_6)$-distinguishing. Considering $\{0, \pm \f 1 6, \pm \f 1 3, \f 1 2\}$ as a copy of $C_6$ sitting inside $S^1$; the only element of $\Aut(C_6)$ which fixes $c_2 :C_6\rightarrow \{R, B\}$ is the reflection about $0$; the only element of $\Aut(C_6)$ which fixes $d_2$ is the reflection about $\f 1 6$; the only element of $\Aut(C_6)$ which fixes $c_3$ is the reflection about $-\f 1 {12}$; the only element of $\Aut(C_6)$ which fixes $d_3$ is the reflection about $\f 1 4$. This means that no two of $\{c_1, c_2, c_3\}$ and $\{d_1, d_2, d_3\}$ can preserve the same reflection, because any such reflection would have to stabilize one element of $W^\prime$, and no $C_6$-reflecition can preserve more than one of the listed colorings.
\\\\ Let $k$ be equal to $c$ or $d$. Assume for the sake of contradiction that none of $k_1, k_2$, or $k_3$ is $O(2)$-distinguishing. 
\\\\ First, suppose that $k_1$ and $k_2$ are invariant under translations $\sigma_1, \sigma_2$. We know that $k_1(\sigma_1(\f 2 3))=R$ while $k_2(\sigma_2(\f 2 3))=B$. This means that $k_2(\sigma_1(\f 2 3))=R$ and $k_1(\sigma_2(\f 1 3))=B$, because Note 2.3.9 tells us that $\sigma_i(W)\cap W = \emptyset$ (as any non-trivial translation which sends elements of $C_6$ to elements of $C_6$ does not even preserve color on $C_6\subset S^1$). Applying $\sigma_2$ and $\sigma_1$, respectively, we see that $k_2(\sigma_2\sigma_1(\f 2 3))=R$ while $k_1(\sigma_1\sigma_2(\f 2 3))=B$. Since $\sigma_1\sigma_2=\sigma_2\sigma_1$, this is only possible if $\sigma_2\sigma_1(\f 2 3)\in W^\prime$. On one hand, $\sigma_2\sigma_1(\f 2 3)\neq \f 2 3$, because then $\sigma_1(\f 2 3)=\sigma_2^{-1}(\f 2 3)$, contradicting the fact that $\sigma_1(\f 2 3)$ and $\sigma_2^{-1}(\f 2 3)$ (which are not elements of $W^\prime$) have different colors under $k$. On the other hand, $\sigma_2\sigma_1(\f 2 3)$ cannot be equal to $0$, because then $\sigma_2\sigma_1(\f 1 6) = \f 1 2$; since $\sigma_1(\f 1 6)\not\in W$, this means that
$$R = k(\f 1 6) = k(\sigma_1(\f 1 6))=k(\sigma_2\sigma_1(\f 1 6))=k(\f 1 2) = B,$$
a contradiction.
\\\\ Next, suppose that $k_1$ is invariant under $\sigma_1$ and $k_2$ is invariant under a reflection, $\tau_2$. Then, we will use the fact that $\sigma_1\tau_2(\f 2 3)= \tau_2\sigma_1^{-1}(\f 2 3)$ to derive a contradiction. Note $2.3.9$ tells us that $\sigma_1(W)\cap W=\emptyset$, and either $\tau_2(0)\not\in W^\prime$ or $\tau_2(\f 2 3)\not\in W^\prime$ (which one depends on whether $k=c$ or $k=d$). Therefore, we know that for some $w\in W^\prime$, 
$$k_2(w)= k_2(\tau_2(w))=k_1(\tau_2(w)) = k_1(\sigma_1\tau_2(w))$$
while
$$k_1(w) = k_1(\sigma_1^{-1}(w))=k_2(\sigma_1^{-1}(w)) = k_2(\tau_2\sigma_1^{-1}(w))$$
implying that $\sigma_1\tau_2(w) \in W^\prime$ (because $k_1(w)\neq k_2(w)$). But $\sigma_1\tau_2(w)\neq w$, because then $\tau_2(w) = \sigma_1^{-1}(w)$, contradicting the fact that these points (which cannot be in $W^\prime$) have different colors under $k$. 
\\\\ Therefore, the only option is that $\sigma_1\tau_2(w)$ is equal to the other element of $W^\prime$. But by replacing $\sigma_1$ with $\sigma_1^2$ (which is guaranteed to be different from $\sigma$ and different from the identity), we can ensure that this does not happen, giving us our contradiction.
\\\\ Since an argument analogous to the above will also work if $k_1$ were invariant under a reflection and $k_2$ were invariant under a translation, we conclude that $k_1$ is preserved by some reflection $\tau_1$ and $k_2$ is preserved by $\tau_2$. Note 2.3.9 tells us that $\tau_1\neq \tau_2$, and thus $\sigma_{21}:= \tau_2\tau_1$ is a nontrivial translation.
\\\\ \textbf{Case 1.} $k=d$. In this situation, we claim that $\sigma_{21}(\f 2 3)\not\in W^\prime$ and $k(\sigma_{21}(\f 2 3))=k_1(\f 2 3) = R$.
\\\\ To prove the claim, note that we know $R = k_1(\f 2 3) = k_1(\tau_1(\f 2 3))=k_2(\tau_1(\f 2 3)) = k_2(\tau_2\tau_1(\f 2 3))$, because $\tau_1(\f 2 3)\not\in W^\prime$ by Note 2.3.9. Therefore, the claim holds provided that $\sigma_{21}(\f 2 3)\not\in W^\prime$. Since $\sigma_{21}$ is nontrivial, we know that $\sigma_{21}(\f 2 3)\neq \f 2 3$. If $\sigma_{21}(\f 2 3) = 0$, then we derive a contradiction from the fact that $\sigma_{21}$ sends the $\f 1 2, \f 5 6, \f 1 6$ triangle to itself. Note 2.3.9 tells us that $\sigma_{21}$ would have to preserve $k$ on this triangle (as no intermediate $\tau_1$-reflection could be an element of $W^\prime$), but this triangle is not monochromatic under $k$. Thus, $\sigma_{21}(\f 2 3)\not \in W^\prime$ and $k(\sigma_{21}(\f 2 3))=R$, as desired.
\\\\ Finally, we consider $k_3$. If $k_3$ is preserved by some translation $\sigma$, then note that $\sigma^{-1}_3(\f 2 3)\not\in W^\prime$ (no $C_6$-translation preserves $k_3$), so
$$B=k_3(\f 2 3) = k_3(\sigma^{-1}_3(\f 2 3)) = k_1(\sigma^{-1}_3(\f 2 3))=k_1(\tau_1\sigma^{-1}_3(\f 2 3))$$
which is equal to $k_2(\tau_1\sigma^{-1}_3(\f 2 3))$ provided that $\tau_1\sigma_3^{-1}(\f 2 3)\not \in W^\prime$. Since $\tau_1(\f 2 3)$ and $\sigma^{-1}_3(\f 2 3)$ have different colors under $k$, we know that $\tau_1\sigma_3^{-1}(\f 2 3)\neq \f 2 3$. By replacing $\sigma_3$ with $\sigma_3^2$ if necessary, we can ensure that $\tau_1\sigma_3^{-1}(\f 2 3)$ is not equal to $0$, which allows us to conclude that $B=k_2(\tau_1\sigma_3^{-1}(\f 2 3))=k_2(\sigma_{21}\sigma_3^{-1}(\f 2 3))$. On the other hand, we know that $\sigma_{21}(\f 2 3)\not\in W^\prime$ and $k(\sigma_{21}(\f 2 3))=R$. This implies that $k_3(\sigma_3^{-1}\sigma_{21}(\f 2 3))=R$, which is a contradiction (as $k_3(\f 2 3)=k_3(0)=B$).
\\\\ If $k_3$ is preserved by a reflection $\tau_3$, then Note 2.3.9 tells us that $\tau_3\neq \tau_1$ and $\tau_3\neq \tau_2$, meaning that $\sigma_{23} := \tau_2\tau_3$ and $\sigma_{31}:= \tau_3\tau_1$ are nontrivial translations, and satisfy the relation
$$\sigma_{21} = \sigma_{23}\sigma_{31}=\sigma_{31}\sigma_{23}.$$
We already know that $\sigma_{21}(\f 2 3)$ is red, and not an element of $W^\prime$. However, $\sigma_{23}(\f 2 3)\neq \f 2 3$ is certainly blue under $k_2$ because $\tau_3(\f 2 3)\not\in W^\prime$ [this can be verified using Note 2.3.9]. Therefore, $\sigma_{23}(\f 2 3)\neq 0$ (which is red under $k_2$), and hence $\sigma_{31}$ sends a blue element of $S^1-W^\prime$ ($\sigma_{23}(\f 2 3)$) to a red element of $S^1-W^\prime$ ($\sigma_{21}(\f 2 3)$); this never happens (it sends the red $\tau_1(\f 2 3)$ to the blue $\tau_3(\f 2 3)$, but never the other way around). Hence, Case 1 leads to a contradition.
\\\\ \textbf{Case 2.} $k=c$. Here, we claim that $\sigma_{21}(0)\not\in W^\prime$ and $k(\sigma_{21}(0)) = B$. To see this, note that we know $B = k_1(0) = k_1(\tau_1(0))= k_2(\tau_1(0))=k_2(\tau_2\tau_1(0))$ [$\tau_1(0)\not\in W^\prime$ by Note 2.3.9], so the claim follows if $\sigma_{21}(0)\not\in W^\prime$. We already know that $\sigma_{21}(0)\neq 0$, and if $\sigma_{21}(0) = \f 2 3$, then $\sigma_{21}$ sends the $\f 1 2, \f 5 6, \f 1 6$ triangle to itself, giving us the same contradiction as in the $k=d$ case.
\\\\ Now, just as in the $k=d$ case, $k_3$ cannot be preserved by a translation $\sigma_3$, because then, after arranging for $\sigma_3\sigma_{21}(0)\not\in W^\prime$, we obtain a contradiction. 
\\\\ If $k_3$ is preserved by a reflection $\tau_3$, then we define $\sigma_{32}$ and $\sigma_{21}$ as before, and use the fact that, $\sigma_{21}(0)$, which is blue under all three extensions of $k$, is also equal to $\sigma_{31}\sigma_{23}(0)$. Furthermore, $\sigma_{23}(0)$ is red under $k_2$, because $\tau_3(0)\not \in W^\prime$ [this can be checked using Note 2.3.9]. Therefore, $\sigma_{23}(0)\neq \f 2 3$ (which is blue under $k_2$), and we conclude that $\sigma_{31}$ sends a red element of $S^1-W^\prime$ to a blue element of $S^1-W^\prime$; this never happens; contradiction.
\\\\ Thus, one of $k_1, k_2$, and $k_3$ is distinguishing with respect to $O(2)$, which proves Lemma 2.3.8. $\hfill \blacksquare$
\\\\ \textbf{Lemma 2.3.10.} \textit{If $W=\{0, \f 1 6, \f 1 3, \f 1 2\}$, then $P(W)$ holds.}
\\\\ \textbf{Proof.} This actually follows from the proof of Lemma 2.3.8 without any extra work. Instead of $W= \{0, \f 1 6, \f 1 3, \f 1 2\}$, we may equivalently consider $W = \{\f 2 3, \f 5 6, 0, \f 1 6\}$. Let $c$ be any precoloring of $S^1-W$ such that $c(\f 1 3)=c(\f 1 2)=B$, and let $d$ be any precoloring of $S^1-W$ such that $c(\f 1 3)=R$ and $c(\f 1 2)=B$. Then, we may extend $c$ and $d$ to $c_1, d_1, c_2, d_2, c_3$, and $d_3$ exactly as in the chart from Lemma 2.3.8. Since we already proved that one of $\{c_1, c_2, c_3\}$ and one of $\{d_1, d_2, d_3\}$ distinguishes $O(2)$, we are done. $\hfill \blacksquare$
\\\\ This completes the proofs of the collection of lemmas necessary for Theorem 2.3.1.
\subsection {Completing the proof of the Theorem 1}
Having finished Theorem 2.3.1, we can loosen the translation constraint further and obtain an even better result for $|W|=4$.
\\\\ \textbf{Condition} $T^{\prime\prime}(W)$: $(W+\f i 5)\cap W = \emptyset$ for $\gcd(i, 5)=1$. 
\\\\ \textbf{Theorem 2.4.1.} \textit{If $W\subset S^1$, $|W|=4$, and $T^{\prime\prime}(W)$ holds, then either $P(W)$ holds, or $W$ is $O(2)$-equivalent to $\{0, \f 1 4, \f 2 4, \f 3 4\}$.}
\\\\ \textbf{Proof.} Suppose $|W|=4$, $T^{\prime\prime}(W)$ holds, and $W$ is not equivalent to $\{0, \f 1 4, \f 2 4, \f 3 4\}$. If $T^\prime(W)$ holds, then by Theorem $2.3.1$, $P(W)$ holds. If $T^\prime(W)$ does not hold, then there are two possibilities.
\\\\ \textbf{Case 1.} $W$ is equivalent to $\{0, \f 1 4, a, b\}$ for $a, b\not\in \{\f 2 4, \f 3 4\}$. In this case, we claim that $R(W)$ holds. In particular, suppose that $R(W)$ does not hold. Then, since $\tau_0(\f 1 4) = -\f 1 4\not\in W$ and $\f 1 2\not\in \{a, b\}$ (i.e., $a$ and $b$ are not fixed by $\tau_0$), we must have that $b = -a$ (or else $\tau_0$ would satisfy the desired property). But since $\tau_{\f 1 4}(0)\not\in W$ and $\f 3 4\not\in \{a, b\}$, we also must have that $b = -a+\f 1 2$; a contradiction. Hence, $R(W)$ holds, and we conclude that $P(W)$ holds by Theorem 2.2.2.
\\\\ \textbf{Case 2.} $W$ is equivalent to $\{0, \f 1 4, \f 3 4, a\}$ for some $a\neq \f 1 2$. Lemma 2.4.2 tells us that $P(W)$ holds in this situation (again, we will present the proof of the lemma below). Assuming this lemma, the proof of Theorem 2.4.1 is complete. $\hfill \blacksquare$
\\\\ \textbf{Lemma 2.4.2.} \textit{If $W \cong \{0, \f 1 4, \f 3 4, a\}$ for $a\neq \f 1 2$, then $P(W)$ holds.}
\\\\ \textbf{Proof.} The principles behind the argument are again similar to those in Lemma 2.3.4, Lemma 2.3.8, and Lemma 2.3.10.  Let $c$ be a precoloring of $S^1-W$; by negating $c$ if necessary, we may assume without loss of generality that $c(\f 1 2) = B$. Then, extend $c$ to $c_1, c_2, c_3$ in the following way: pick $c_1(a)=c_2(a)=c_3(a)$ such that exactly three of $\{a, -a, a+\f 1 2, -a+\f 1 2\}$ have the same color (this is always possible because either $-a, a+\f 1 2$, and $-a+\f 1 2$ all have the same color to begin with, or exactly two of the three have the same color). This choice guarantees that $c_1, c_2$, and $c_3$ are not preserved by $\sigma_{\f 1 2}$, $\tau_0$, or $\tau_{\f 1 4}$. The rest of $W$ is colored according to this chart.
\begin{center}
  \begin{tabular}{ | l || c | c | c |}
    \hline
    & $c_1$ & $c_3$ & $c_2$ \\ \hline
    $0$ & R & B & B \\ \hline
    $\f 1 4$ & B & B & R \\ \hline
    $\f 3 4$ & R & R & R \\ \hline
  \end{tabular}
\end{center}
Note that $c_1$ and $c_2$ only differ on $W^\prime := \{0, \f 1 4\}$. Also note that $\sigma_{\f 1 2}$, $\sigma_{\f 1 4}$, $\tau_0$, $\tau_{\f 1 4}$, and $\tau_{\f 1 8}$ do not preserve any of the $c_i$.. We will now run through the same argument as before -- assuming that $c_1, c_2$, and $c_3$ all satisfy some symmetry, showing that all of the symmetries must be reflections $\tau_1, \tau_2, \tau_3$, and deriving a contradiction from the fact that $\tau_2\tau_1 = (\tau_2\tau_3)(\tau_3\tau_1)=(\tau_3\tau_1)(\tau_2\tau_3)$.
\\\\ Assume for the sake of contradiction that $P(W)$ does not hold; in particular, we assume that none of $c_1$, $c_2$, and $c_3$ are distinguishing. Suppose that $c_1$ and $c_2$ are both preserved by nontrivial translations $\sigma_1, \sigma_2$. Then, $c_1(\sigma_1(0)) = c_1(0)=R$, which implies that $\sigma_1(0)\neq \f 1 4$. Thus, $\sigma_1(0)\not\in W^\prime$, and hence $c_2(\sigma_2\sigma_1(0))=c_2(\sigma_1(0))=c_1(\sigma_1(0))=R$. This argument is symmetric in $c_1, c_2$, so we also conclude that $c_1(\sigma_2\sigma_1(0))=B$. This can only happen if $\sigma_2\sigma_1(0))=\f 1 4$, but this implies that $\sigma_2\sigma_1(\f 1 2) = \f 3 4$. Since $\sigma_1(\f 1 2)\not\in W^\prime$ (as $\sigma_{\f 1 2}$ and $\sigma_{\f 1 4}$ do not preserve $c_1$), we know that $c(\sigma_1(\f 1 2))=c_1(\f 1 2)=B$. Then, we have $\sigma_2$ taking a blue element of $S^1-W^\prime$ to a red element of $S^1-W^\prime$; a contradiction. Therefore, it is not the case that $c_1$ and $c_2$ are both preserved by translations.
\\\\ Now suppose that $c_1$ is preserved by a translation $\sigma_1$ and $c_2$ is preserved by a reflection $\tau_2$. Since $\tau_2(0)\not \in W^\prime$ ($\tau_0$ and $\tau_{\f 1 8}$ do not preserve $c_1$ or $c_2$), we know that $c(\tau_2(0))=c_2(0)=B$, and hence $c_1(\sigma_1\tau_2(0))=B$. Similarly, since $\sigma^{-1}_1(0)\not\in W^\prime$, we have that $c_2(\tau_2\sigma_1^{-1}(0))=c(\sigma_1^{-1}(0))=c_1(0)=R$. Since $\tau_2\sigma^{-1}_1=\sigma_1\tau_2$, this is only possible if $\tau_2\sigma_1^{-1}(0)=\f 1 4$. But then we would have that $\tau_2\sigma_1^{-1}(\f 1 2)=\f 3 4$, which is a contradiction (lack of color preservation) because $\sigma_1^{-1}(\f 1 2)\not\in W^\prime$. The same argument applies when $c_1$ is preserved by a reflection and $c_2$ is preserved by a translation. 
\\\\ We conclude that $c_1$ and $c_2$ are both preserved by reflections $\tau_1$ and $\tau_2$. Since $\sigma_{\f 1 2}$ and $\sigma_{\f 1 4}$ also do not preserve $c_3$, the argument from the previous paragraph also proves that $c_3$ must also be preserved by a reflection, $\tau_3$. Because none of the colorings are preserved by $\tau_0$ or $\tau_{\f 1 4}$, we know that $\tau_1, \tau_2$, and $\tau_3$ are pairwise distinct, allowing us to define nontrivial translations $\sigma_{ij} := \tau_i\tau_j$.
\\\\ We already know that $\tau_i(W^\prime)\cap W^\prime = \emptyset$ for all $i$, so we have that $c_i(\sigma_{ij}(0))=c(\tau_j(0))=c_j(0)$. Furthermore, we also know that $\sigma_{ij}(0)\not\in W^\prime$, because $\sigma_{ij}(0)=\f 1 4\rightarrow \sigma_{ij}(\f 1 2) = \f 3 4$, which is only possible if $\tau_j(\f 1 2)\in W^\prime$. But $\tau_j(\f 1 2)\neq 0$ (because $\tau_{\f 1 4}$ never preserves a coloring) and $\tau_j(\f 1 2)=\f 1 4\rightarrow \tau_i(\f 1 4) = \f 3 4$, a contradiction ($\tau_0$ never preserves a coloring). Thus, $c(\sigma_{ij}(0))=c_j(0)$ for all $i\neq j$. In particular, $c(\sigma_{23}(0))=B$ and $c(\sigma_{21}(0))=R$. But $\sigma_{31}\sigma_{23}(0)=\sigma_{21}(0)$, so $\sigma_{31}$ sends a blue element of $S^1-W^\prime$ to a red element of $S^1-W^\prime$. Since $c_3$ and $c_1$ differ only at $0$ (where $c_1(0)=R$ and $c_3(0)=B$), it is easy to see that this is impossible. Thus, one of $c_1, c_2$, and $c_3$ is distinguishing, which proves Lemma 2.4.2. $\hfill \blacksquare$ 
\\\\ \textbf{Corollary 2.4.3.} \textit{If $W\subset S^1$, $|W|=5$, and $T^{\prime\prime}(W)$ holds, then $P(W)$ holds.}
\\\\ \textbf{Proof.} If $W\subset S^1$, $|W|=5$, and $T^{\prime\prime}(W)$ holds, then there is some subset $W_4$ of $W$ of size $4$ which is not equivalent to $\{0, \f 1 4, \f 1 2, \f 3 4\}$ (at worst, $W$ can be equivalent to $\{0, \f 1 4, \f 1 2, \f 3 4, a\}$, and we can remove $0$, for instance). Then, by Theorem 2.4.1, $P(W_4)$ holds, and hence $P(W)$ holds. $\hfill \blacksquare$
\\\\ Finally, we will remove all constraints on $W$ to prove Theorem 1.
\\\\ \textbf{Theorem 2.4.4.} \textit{If $W\subset S^1$ and $|W|=5$, then $P(W)$ holds unless $W\cong \{0, \f 1 5, \f 2 5, \f 3 5, \f 4 5\}$.}
\\\\ \textbf{Corollary 2.4.5.} \textit{If $W\subset S^1$ and $|W|=6$, then $P(W)$ holds unconditionally.}
\\\\ \textbf{Proof.} Let $W\subset S^1$ be such that $|W|=5$ and $W$ is not equivalent to the one exceptional set, and let $W_4$ be equal to \textit{any} subset of $W$ of size $4$. Then, if $T^{\prime\prime}(W_4)$ holds, either $P(W_4)$ holds (in which case $P(W)$ holds), or $W_4\cong \{0, \f 1 4, \f 2 4, \f 3 4\}$. In the latter case, let $a\in W-W_4$; then, $P(\{a\}\cup W_4- \{\f 2 4\})$ holds by Lemma 2.4.2. Therefore, $T^{\prime\prime}(W_4)$ cannot hold for any such $W_4$. Then, there are two possibilities.
\\\\ \textbf{Case 1.} $W\cong \{0, \f 1 5, \pm \f 2 5, a, b\}$ for some $a, b\in S^1$ [the $\pm $ notation here means that \textit{either} $\f 2 5\in W$ or $-\f 2 5\in W$, but not both]. In this case, without loss of generality, let $a$ be such that $a\not\in \{\pm \f 3 5, \f 4 5\}$ (either $a$ or $b$ must satisfy this property). Let $W_4 = \{0, \f 1 5, \pm \f 2 5, a\}$. We claim that $R(W_4)$ holds. To verify this, suppose that it were not the case. Since $\tau_0(\f 1 5)=-\f 1 5$ and $\tau_0(\pm\f 2 5) = \mp \f 2 5$ are not in $W_4$, if $R(W_4)$ does not hold, then $-a\in W_4$. But we know that $-a\not\in \{0, \f 1 5, \pm \f 2 5\}$, so we conclude that $-a = a$ and hence $a=\f 1 2$. But then we can consider $\tau_{\f 1 5}$ and $\tau_{\pm \f 2 5}$; one of these reflections ($\tau_{\f 2 5}$ in the $\f 2 5$ case, and $\tau_{\f 1 5}$ in the $-\f 2 5$ case) satisfies $\tau(W_4-\{a\})\cap W_4=\emptyset$. Therefore, if $R(W_4)$ does not hold, then we would need for $\tau(a)=a$ for this second reflection as well; since $\f 1 2$ is not invariant under either $\tau_{\f 1 5}$ or $\tau_{\f 2 5}$, we conclude that $R(W_4)$ holds. As a result, Theorem 2.2.2 tells us that $P(W_4)$ holds, and thus $P(W)$ holds.
\\\\ \textbf{Case 2.} $W\cong \{0, \f {i_0}5, a, a+\f {i_a}5, b\}$ for some $a, b\in S^1$ such that the $\sigma_{\f 1 5}$-orbits of $0, a$, and $b$ are all distinct. In this case, let $W_3 = \{0, a, b\}$; if $\tau_0(W_3-\{0\})\cap W_3 = \emptyset$, then we define $W_4 = W_3 \cup \{\f {i_0}5\}$. Since $\tau_{w_0}(\f{i_0}5) = - \f{i_{w_0}}5$ (which is not an element of $W$ by the $\sigma_{\f 1 5}$-orbit condition), we conclude that $R(W_4)$ holds, and hence $P(W_4)$ holds by Theorem 2.2.2. The same reasoning applies if $\tau_a(W_3 - \{a\})\cap W_3=\emptyset$. If neither the $\tau_0$ condition nor the $\tau_a$ condition holds, it is easily verified that $a = \f 1 2$ or $W_3$ is equivalent to either $\{0, \f 1 4, \f 3 4\}$ or $\{0, \f 1 3, \f 2 3\}$. If $a=\f 1 2$, then we let $W_4 = \{0, \f {i_0}5, \f 1 2 + \f{i_{\f 1 2}}5, b\}$ and note that $R(W_4)$ holds (in particular, $\tau_0(W_4-\{0\})\cap W_4 =\emptyset$ by the $\sigma_{\f 1 5}$-condition). Hence, $P(W_4)$ holds by Theorem 2.2.2, and so $P(W)$ holds. If $W_3 \cong \{0, \f 1 4, \f 3 4\}$, then $P(W_3\cup \{\f {i_0}5\})$ holds by Lemma 2.4.2, and hence $P(W)$ holds. Finally, we can ensure that $W_3\neq \{0, \f 1 3, \f 2 3\}$ by applying a translation by $\{-\f {i_0}5\}$ if necessary (i.e., considering $\{\f {i_0}5, a, b\}$ if $\{0, a, b\} = \{0, \f 1 3, \f 2 3\}$). Therefore, we may assume without loss of generality that $R(W_3)$ holds (in all cases in which we are not already done). Let $w_0\in W_3$ be as described in condition $R(W_3)$, and let $W_4 = W_3 \cup \{w_0 + \f{i_{w_0}}5\}\subset W$.  Thus, since $R(W_4)$ holds, $P(W_4)$ holds by Theorem 2.2.2, and hence $P(W)$ holds in all cases. $\hfill \blacksquare$
\subsection{An interesting corollary}
As a result of our work above (in particular, from our proof of Lemma 2.3.4), we also get a result of a slightly different flavor.
\\\\ \textbf{Corollary 2.5.1.} Let $c$ be any $2$-coloring of $S^1$. Then, there exists a distinguishing coloring $c^*: S^1\rightarrow \{R, B\}$ such that $c(x)\neq c^*(x)$ for at most three values of $x$. 
\\\\ \textbf{Proof.} Suppose we have a $2$-coloring of $S^1$, denoted $c$. If $c$ is identically one color, then changing the colors of $0, \f 1 3$, and $\f 1 2$ suffices to produce a distinguishing $2$-coloring. 
\\\\ If $c$ is not identically one color, then we claim that there exists a reflection (about some point $w_0$) $\tau_{w_0}$ and a point $a\neq w_0 \pm \f 1 4$ such that $c(\tau_{w_0}(a))\neq c(a)$. If this were not the case, then, for example, $c(0)$ must be equal to $c(\theta)$ for all $\theta \neq \pm \f 1 4$. Similarly, $c(\f 1 3)$ must be equal to $c(\theta)$ for all $\theta\neq \f 1 3 \pm \f 1 4$. But then by transitivity we find that $c(0)=c(\theta)$ for all $\theta\in S^1$, contradicting the fact that $c$ is not uniformly one color. This proves the claim.
\\\\ Let $\tau_{w_0}$ be a reflection described in the claim. By translational symmetry, we may assume that $w_0=0$, so there exists some $a\neq \pm \f 1 4$ such that $c(a)\neq c(-a)$. Let $W = \{0, a, \f 1 2, a+\f 1 2\}$, and let $c^\prime$ be the restriction of $c$ to $S^1-W$. Lemma $2.3.4$ tells us that there exists an extension $c^*$ of $c^\prime$ which is distinguishing; furthermore, the proof of Lemma 2.3.4 specifies that there exists a distinguishing $c^*$ such that $c^*(a)\neq c^\prime(-a)= c(-a)$. This means that $c^*(a)=c(a)$, so $c^*$ and $c$ differ at most on $\{0, \f 1 2, a+\f 1 2\}$. $\hfill \blacksquare$
\section{Extending precolorings on $\R^2$: a proof of Theorem 3}
The complexity of extending precolorings on $\R^2$ is highly dependent on the choice of symmetry group $\Gamma$. First, we will show that the case of $\Gamma = O(2)$ has already been resolved by Theorem 1.
\\\\ \textbf{Theorem 3.1.} \textit{$\ext_D(\R^2, O(2))=7$.}
\\\\ \textbf{Proof.} The fact that $7$ is a lower bound to the extension number was proved in [2]. Let $W\subset \R^2$ be such that $|W|=7$ and the pointwise stabilizer $\stab_{O(2)}(W)$ is trivial, and let $c$ be a precoloring of $\R^2-W$. Assume for the sake of contradiction that this precoloring cannot be extended to a distinguishing coloring of $\R^2$. Note that the action of $O(2)$ on $\R^2 = \bigcup_{r\in \R_{\geq 0}} r\cdot S^1$ can be decomposed into separate actions of $O(2)$ on each individual $r\cdot S^1$. If there is any $r\in \R$ such that $|W\cap r\cdot S^1|\geq 6$, then $P(W, O(2))$ holds by Theorem 1. If not, then there are at least two nondegenerate circles in $\R^2$ which intersect $W$.
\\\\ We now claim that there exist two points $x_1, x_2\in W-\{0\}$ such that $x_1$ and $x_2$ are not on the same $r\cdot S^1$ and the line connecting $x_1$ to $x_2$ does not intersect $0\in \R^2$. Suppose this were not the case. Let $C=r\cdot S^1$ be a nondegenerate circle such that $|C\cap W|>0$ is minimal among circles that intersect $W$. If $C\cap W = \{x_1\}$, then if the claim is false, all of the other points in $W$ lie on the line connecting $0$ to $x_1$; this contradicts the stabilizer condition, because reflection across this line pointwise stabilizes $W$. On the other hand, if $C\cap W\supset \{x_1, x_2\}$ and the claim is false, we get that all elements of $W-C\cap W$ lie on the line connecting $0$ to $x_1$ as well as the line connecting $0$ to $x_2$, showing that $x_2$ also lies on the line connecting $0$ to $x_1$. Applying this reasoning to all pairs of points inside $C\cap W$, we again obtain that $W$ lies on a line, a contradiction. Thus, we may find $x_1, x_2$ as stated.
\\\\ Let $x_1\in r_1\cdot S^1$ and $x_2\in r_2\cdot S^1$ be two points which satisfy the claim. Color the rest of $W$ red (or any other combination of colors). Then, let $c_1$ be the coloring of $r_1\cdot S^1$ where $x_1$ is red, and $c_2$ be the coloring where $x_1$ is blue. If $c_1$ satisfies an $SO(2)$ symmetry $\sigma_1$ and $c_2$ satisfies any $O(2)$ symmetry which does not fix $x_1$, ($\sigma_2$ or $\tau_2$) we obtain the usual contradiction from the relation $\sigma_1\sigma_2=\sigma_2\sigma_1$ or $\sigma_1\tau_2 = \tau_2\sigma_1^{-1}$. The same holds where $c_1$ and $c_2$ are exchanged. Therefore, we may assume that $c_1$ satisfies a reflection $\tau_1$ which fixes $x_1$ and no other symetries, or both $c_1$ and $c_2$ satisfy reflections $\tau_1, \tau_2$ and no other symmetries.
\\\\ In either case, $c_1$ satisfies only $\tau_1$, a reflection. If $\tau_1(x_2)\neq x_2$, we may color $x_2$ so that the final coloring $c^*$ distinguishes $O(2)$, a contradiction. If $\tau_1(x_1)=x_1$, then by the claim, $\tau_1(x_2)\neq x_2$, so we are done. On the other hand, if $\tau_1(x_1)\neq x_1$, then $c_2$ satisfies only $\tau_2\neq \tau_1$, so one of $\tau_1$ or $\tau_2$ must satisfy $\tau_i(x_2)\neq x_2$. Thus, we have proved Theorem 3.1. $\hfill \blacksquare$
\\\\ The full isometry group $E(2)\supset O(2)$ is much more difficult to deal with. First, we'll classify the elements of $E(2)$, based on the identification $\R^2\cong \C$: every $\gamma\in E(2)$ is either a translation ($z\in \C\mapsto z+a$), a rotation about some point ($z\mapsto \omega z + a$, $\omega\in S^1, a\in \C$), a reflection over some line ($z\mapsto \f{-a}{\overline a} \overline z+a, a\in \C$), or a ``glide reflection'', which is a reflection over a line combined with a translation parallel to that line. The subgroup of $E(2)$ composed of translations and rotations is called $SE(2)$. 
\\\\ So far, we have used the work done on $S^1$ to prove $\ext_D(\R^2, O(2))=7$. We can also apply the work from [2] done on $(\R, E(1))$ by considering the following subgroup $\Gamma$ of $E(2)$: the group of translations $z\mapsto z+a$ and $180^\circ$ rotations $z\mapsto -z+2a$. By using the techniques from [2], we can prove the following lemma.
\\\\ \textbf{Lemma 3.2.} \textit{Let $\Gamma\subset E(2)$ be as described above. Then, $\ext_D(\R^2, \Gamma)=4$}.
\\\\ \textbf{Proof.} The fact that $4$ is a lower bound to the extension number follows from the fact that $\ext_D(\R, E(1))=4$; whenever $W$ is contained in a line $\ell$, a precoloring of $\R^2 - W$ which colors $\R^2-\ell$ red must be extended to distinguish the action of $E(1)$ on that line.
\\\\ Let $W\subset \R^2$ be of size 4, and let $c$ be a precoloring of $\R^2-W$. Assume that $c$ cannot be extended to a distinguishing coloring of $\Gamma$. Recall the definition of the property $R(W)$ as it pertains to $\Gamma$: $R(W)$ holds if there exists some $w_0\in W$ such that $\tau_{w_0}$ (which is now the $180^\circ$ rotation about $w_0$) sends $W-\{w_0\}$ outside of $W$. 
\\\\ \textbf{Claim.} $R(W)$ always holds. 
\\\\ \textbf{Proof.} Among all $w\in W$ with minimal $x$-coordinate, pick the unique $w_0$ with minimal $y$-coordinate. Drawing axes perpendicular to the $x$ and $y$-axes which meet at $w_0$, it is clear that all of $W$ sits inside the union of the right half-plane and the positive $y$-axis as defined by those axes. Therefore, $\tau_{w_0}(W)$ sits inside the left half-plane and the negative $y$-axis; this means that $\tau_{w_0}(W)\cap W \subset \{w_0\}$, proving the claim.
\\\\ The rest of the proof of Lemma 3.2 follows exactly as the proof of Theorem 2.2.2, with no modifications. $\hfill \blacksquare$
\\\\ We will now use Lemma 3.2 to prove the second half of Theorem 3. Let $SE(2)= \{\vec v \mapsto A \vec v + \vec b, A\in SO(2), \vec b\in \R^2\}$.
\\\\ \textbf{Theorem 3.3.} \textit{$\ext_D(\R^2, SE(2))<\infty$.}
\\\\ Before proving the theorem, we first analyze the case of $|W|=4$ in detail. 
\\\\ \textbf{Lemma 3.4.}  Let $W\subset \R^2$ be such that $|W|=4$, and let $c$ be a precoloring of $\R^2-W$. Suppose that $c$ cannot be extended to a $SE(2)$-distinguishing coloring of $\R^2$. Then, there exists an extension $c_0$ of $c$ which is $\Gamma$-distinguishing, and satisfies a $120^\circ$ rotational symmetry.
\\\\ \textbf{Proof.} Let $W$ and $c$ be as described in the statement of the lemma. The following technical lemma will be used to derive a contradiction.
\\\\ \textbf{Lemma 3.5.} There exists an extension $c_0$ of $c$ which is $\Gamma$-distinguishing, and preserved by some rotation $\gamma_1$ of either odd or infinite order (without loss of generality, about the origin $0\in \R^2$). Furthermore, for at least one such $c_0$, there exists an extension $c_1$ of $c$ which is preserved under $\gamma_1$ and a point $x_0\in W$ which is not fixed by $\gamma_1$, such that the extension $c_2$ obtained from switching $c_1(x_0)$ to the opposite color is not preserved by either a rotation about $x_0$ or the map $z\mapsto -z+x_0$. 
\\\\ \textbf{Proof:} First, note that if we choose $c_1$ so that it is $\Gamma$-distinguishing and rotationally symmetric about $0$, then changing $c_1(x_0)$ to the opposite color for any $x_0\in W-\{0\}$ will never result in a rotational symmetry about $x_0$. This is because if such a symmetry did result, then $c_1$ would itself be symmetric under some rotation about $x_0$ as well as another rotation about $0$. The commutator of these two rotations is a nontrivial translation, so $c_1$ would not be $E_1$-distinguishing. Therefore, we only have to deal with $x_0$-rotational symmetry if $c_1$ is not chosen to be $\Gamma$-distinguishing. 
\\\\ Let $c_0$ be an extension of $c$ to $\R^2$ which is $\Gamma$-distinguishing (Lemma 3.2 guarantees that $c_0$ exists). Since $c_0$ cannot be $SE(2)$-distinguishing, $c_0$ must be preserved by some rotation which is not $180^\circ$. Without loss of generality, we may assume that the rotation is about the point $0\in \R^2$, so $c_0$ is preserved by $\gamma_0: z\mapsto \omega z$. Furthermore, $\gamma_0$ cannot have even order, for otherwise some power of $\gamma_0$ is the map $z\mapsto -z$, which we know does not preserve $c_0$. 
\\\\ Now, suppose that for every $x\in W-\{0\}$, the coloring $c_x$ obtained from changing $c_0(x)$ to the opposite color \textit{is} preserved by $z\mapsto -z+x$. 
\\\\ \textbf{Case 1.} $0\not\in W$. Then for every $x\in W$, we have that $c_x(x) = c_x(0)=c_0(0)$, and so $c_0(x)\neq c_0(0)$. In other words, all of the points in $W$ have the same color under $c_0$. Now, let $x_1, x_2, x_3\in W$ be such that $x_1\neq 2^{\pm 1} x_2$ and $x_1\neq 2^{\pm 1} x_3$ (three such points in $W$ exist). Let $c_1$ be the coloring which matches $c_0$ except at $x_1$ and $x_2$. Furthermore, by switching $x_2$ and $x_3$ if necessary (in the definition of $c_1$), we may assume that the $180^\circ$ rotation about $0$ does not preserve $c_1$. 
\\\\ We want to show that $c_1$ satisfies the properties demanded by Lemma 3.5. In particular, we claim that $c_1$ is $\Gamma$-distinguishing. 
\\\\ To show this, suppose that $c_1$ is preserved by a translation $\sigma_1: z\mapsto z+a$. By replacing $a$ with a sufficiently large multiple of $a$, we can ensure that $a+x_1-x_2\not\in \{x_1, x_2\}$ and $x_2-a\neq x_1$ (so it makes sense to talk about $c(a+x_1-x_2)$ and $c(x_2-a)$). Furthermore, we know already that $x_1-x_2\not\in \{x_1, x_2\}$ by our choice of $x_1$ and $x_2$. Therefore, we see that $c(a+x_1-x_2) = c(\sigma_1(x_1-x_2))=c(x_1-x_2)=c_{x_1}(x_2)=c_0(x_2)$, but also that $c(a+x_1-x_2) = c(-(x_2-a)+x_1) = c(x_2-a) = c_1(x_2)\neq c_0(x_2)$, a contradiction. 
\\\\ On the other hand, if $c_1$ is preserved by a $180^\circ$ rotation $\gamma_1: z\mapsto -z+a$, by the construction of $c_1$ we know that $a\not\in \{0, x_1, x_2\}$. Furthermore, we know that $x_1-x_2\not\in \{x_1, x_2\}$ Therefore, we see that $c_{x_1}(x_1-a)=c(a) = c(0)=c_1(x_2)\neq c_{x_1}(x_2)=c(x_1-x_2)=c_1(x_2-x_1+a)$. Now, $c_1(x_2-x_1+a)=c_{x_2}(x_2-x_1+a)$ as long as $x_2-x_1+a\neq x_1$, and we already know that $c_1(x_1)=c_1(x_2)\neq c_{x_1}(x_2)$, so this cannot happen. Therefore, we have that $c_{x_1}(x_1-a)=c(0)\neq c_1(x_2-x_1+a)= c_{x_2}(x_2-x_1+a)=c_{x_2}(x_1-a)$, which can only happen if $x_1-a=x_1$, i.e., $a=0$; a contradiction. This proves the claim that $c_1$ is $\Gamma$-distinguishing. 
\\\\ Since $c_1$ is $\Gamma$-distinguishing, we obtain $\gamma_1$ analogously to $\gamma_0$ as described before. Since two elements of $W$ are colored red under $c_1$ and two elements are colored blue under $c_1$, it is now certain that we can pick $x_0\in W-\stab_{\gamma_1}(\R^2)$ as the lemma describes.
\\\\ \textbf{Case 2.} $0\in W$. First of all, the argument from Case 1 still works almost all of the time. In this situation, we have that the three points in $W-\{0\}$ are all the same color (and opposite the color of $0$) under $c_0$. We can still flip the colors of $x_1, x_2\in W-\{0\}$ to obtain another $\Gamma$-distinguishing coloring unless the following things both happen.
\\\\ 1) $W = \{0, \f 1 2 x, x ,2x\}$ for some $x$ (meaning we can only find two points $x_1, x_2$, and not a third).
\\ 2) $c_0(-\f 1 2 x) = c_0(-2 x)=c_0(0)$ and $c_0(-x)=c_0(x)$ (when we change the colors of $\f 1 2 x_0$ and $2x_0$, a $z\mapsto -z$ symmetry results).
\\\\ However, if (1) and (2) both happen, we note that $c_0(-\f 1 2 x)\neq c_0(x)$, contradicting the fact that flipping the color of $\f 1 2 x$ is supposed to result in a $z\mapsto -z+\f 1 2 x$ symmetry. Therefore, we can still find $x_1, x_2\in W-\{0\}$ such that flipping the $c_0$-colors of $x_1$ and $x_2$ results in an $\Gamma$-distinguishing coloring. Under this new coloring, which we call $c^*$, $0, x_1$, and $x_2$ all have the same color while $x_3$ (the last element of $W$) has the opposite color. As long as $c^*$ is not rotationally symmetric about $x_3$, we are again done.
\\\\ Therefore, we may assume that $c^*$ is symmetric under some rotation about $x_3$. Then, pick $c_1$ to be the coloring which matches $c_0$ except at $0$ (this is still symmetric under $\gamma_0$) and pick $x_0 = x_3$. Since $c_1$ itself is symmetric under $z\mapsto -z+x_3$, we know that $c_2$ (which is equal to $c_1$ except at $x_3$) is not symmetric under this rotation. Then, we are done unless $c_2$ is symmetric under some rotation about $x_3$. This, combined with the fact that $c^*$ is also symmetric under some rotation about $x_3$, gives us a contradiction (from the usual commutation relation) unless $0, x_1$, and $x_2$ lie on a circle centered at $x_3$ such that $0, x_1$, and $x_2$ form an equilateral triangle.
\\\\ Finally, this condition is actually symmetric under exchange of $c_0$ and $c^*$ as the intial coloring (because both are $\Gamma$-distinguishing), so $x_3$ would also have to form an equilateral triangle with two out of $\{0, x_1, x_2\}$. However, the first triangle being equilateral and centered around $x_3$ rules out the possibility that any such second triangle could be equilateral; contradiction. This proves Lemma 3.5. $\hfill \blacksquare$
\\\\ To finish the proof of Lemma 3.4, let $c_1$, $\gamma_1: z\mapsto \omega_1 z$, $x_0$, and $c_2$ be as asserted by Lemma 3.5. Since we assumed that $c$ is not $SE(2)$-distinguishing, $c_2$ must be preserved by some $\gamma_2: z\mapsto \omega_2 z + a$, where $\gamma_2(x_0)\neq x_0$ and $(\omega_2, a)\neq (-1, x_0)$ by the lemma. If $a=0$, since $c_1$ and $c_2$ differ only at $x_0$, we can derive the usual contradiction from the fact that $\gamma_1\gamma_2=\gamma_2\gamma_1$; therefore, we may suppose that $a\neq 0$. Now, note that the commutator $\sigma_1 = \gamma_2^{-1}\gamma_1^{-1}\gamma_2\gamma_1$ is a translation, $z\mapsto z+\f{a(1-\om_1)}{\om_1\om_2}$. Since we know that $a\neq 0$ and $\om_1\neq 1$, we in fact know that this is a nontrivial translation. Similarly, we have that $\sigma_2 := \gamma_1^{-2}\gamma_2^{-1}\gamma_1^2\gamma_2: z\mapsto z - \f{a(1-\om_1^2)}{\om_1^2\om_2}$ is a nontrivial translation, as $a\neq 0$ and $\om_1^2\neq 1$. Therefore, we have the relation $\sigma_1\sigma_2=\sigma_2\sigma_1$, which in terms of $\gamma_1$ and $\gamma_2$ becomes
$$\gamma_2^{-1}\gamma_1^{-1}\gamma_2\gamma_1^{-1}\gamma_2^{-1}\gamma_1^2\gamma_2 =\gamma_1^{-2}\gamma_2^{-1}\gamma_1\gamma_2\gamma_1.$$
Without loss of generality, suppose that $c_1(x_0)=R$. Noting that $\sigma_1\sigma_2(x_0)=x_0+\f{a(\om_1-1)}{\om_1^2\om_2}\neq x_0$, we will show that, under the caveat that $\gamma_1$ may be replaced by $\gamma_1^{n}$ for some $n$ and $\gamma_2$ may be replaced with $\gamma_2^2$, $c(\sigma_1\sigma_2(x_0))=B$ while $c(\sigma_2\sigma_1(x_0))=R$.
\\\\ We start by looking at the right side of the equation (and applying each letter of the word one at a time). We know that $\gamma_1(x_0)=\om_1x_0\neq x_0$, and hence $c(\gamma_1(x_0))=R$. Therefore, we have that $c(\sigma_2\sigma_1(x_0))=R$ unless one of the following two things happens:
\\\\ 1) $\gamma_1\gamma_2\gamma_1(x_0)=x_0$.
\\ 2) $\sigma_2\sigma_1(x_0)=x_0$ (which we have already ruled out).
\\\\ It does not matter if, say, $\gamma_2\gamma_1(x_0)=x_0$, because then $\gamma_2\gamma_1(x_0)$ will still be red under $c_1$, and we can continue applying the next letter. Writing out an explicit formula for $\gamma_1\gamma_2\gamma_1$ then tells us that $c(\sigma_2\sigma_1(x_0))=R$ unless $\om_1^2\om_2x_0+a\om_1 = x_0$, i.e., $x_0 = \f{a\om_1}{1-\om_1^2\om_2}:= A_1(\gamma_1, \gamma_2)$.
\\\\ Now, analyzing the left hand side of the equation, we know that $\gamma_2(x_0)\neq x_0$ by Lemma 3.5, so $c(\gamma_2(x_0))=B$. Therefore, we have that $c(\sigma_1\sigma_2(x_0))=B$ unless one of the following two things happens:
\\\\ 1) $\gamma_2^{-1}\gamma_1^2\gamma_2(x_0)=x_0$, which reduces to the equation $x_0 = \f{-a}{\om_2} := A_3(\gamma_1, \gamma_2)$.
\\ 2) $\gamma_2\gamma_1^{-1}\gamma_2^{-1}\gamma_1^2\gamma_2(x_0)=x_0$, which reduces to the equation $x_0 = \f{a(\om_1^2+\om_1-1)}{\om_1(1-\om_1\om_2)}:= A_2(\gamma_1, \gamma_2)$.
\\\\ Thus, we obtain our contradiction unless $x_0 = A_i(\gamma_1, \gamma_2)$ for some $1\leq i\leq 3$. Now, we show that with the proper modifications, we can ensure that this does not happen.
\\\\ \textbf{Step 1.} Possibly replacing $\gamma_2$ with $\gamma_2^2$, we can ensure that $x_0\neq A_3(\gamma_1, \gamma_2)$.
\\\\ Suppose that $x_0 = \f{-a}{\om_2}$. Then, as $\gamma_2^2(z) = \om_2^2 z + \om_2 a + a$, we see that $A_3(\gamma_1, \gamma_2^2) = \f{-a(\om_2+1)}{\om_2^2}$, and if $x_0 = A_3(\gamma_1, \gamma_2^2)$ as well, we obtain that either $a=0$ or $\om_2 = \om_2+1$; a contradiction in either case. Therefore, as long as $\gamma_2^2$ does not fix $x_0$, we can safely replace $\gamma_2$ with $\gamma_2^2$. Furthermore, since $\gamma_2^2$ is a rotation about some $y\neq x_0$, $\gamma_2^2$ fixes $x_0$ if and only if either $\om_2=-1$ or $\gamma_2^2$ is the identity (which also implies that $\om_2 = -1$). In this situation, since $x_0 = \f{-a}{\om_2}$, we get that $x_0 = a$, which (combined with $\om_2 = -1$) cannot be the case by Lemma 3.5. Thus, we may assume without loss of generality that $x_0\neq A_3(\gamma_1, \gamma_2)$. Furthermore, since $A_3(\gamma_1, \gamma_2)$ is actually independent of $\gamma_1$, any future modifications to $\gamma_1$ will not change this fact.
\\\\ \textbf{Step 2.} Attempt to replace $\gamma_1$ with $\gamma_1^2$.
\\\\ Since $\gamma_1$ does not have even order, $\gamma_1^2$ also does not have even order, so we may repeat the above process using $\gamma_1^2$ instead of $\gamma_1$. Therefore, we obtain our desired contradiction unless $x_0 = A_i(\gamma_1, \gamma_2)$ for some $1\leq i\leq 2$ and $x_0 = A_j(\gamma_1^2, \gamma_2)$ for some $1\leq j\leq 2$. This gives us four pairs of simultaneous equations in the three variables $\om_1, \om_2, a$, whose solutions (simplified using $a\neq 0$ and $\om_1 \not\in\{0, 1\}$) are as follows:
$$B_{1, 1}(\om_1, \om_2, a) := \om_2 + \f 1 {\om_1^3} = 0,$$ 
$$B_{1, 2} (\om_1, \om_2, a) := \om_1^3 + \om_1 + 1 = 0,$$
$$ B_{2, 1}(\om_1, \om_2, a) := \om_2 + \f{\om_1^2 - 1}{\om_1^4 (\om_1+2)} = 0,$$
$$B_{2, 2}(\om_1, \om_2, a) := \om_2 + \f{\om_1^3 + 1}{\om_1(\om_1^2 - \om_1 - 1)} = 0.$$
First, note that for a fixed $\om_2$, there are at most $3+3+5+3 = 14$ values of $\om_1$ which satisfy any of the four above equations. Leaving $\om_2$ fixed, we may now replace $\om_1$ (respectively, $\gamma_1$) with $\om_1^n$ ($\gamma_1^n$) for any $n\in \Z$ unless $\om_1^{2n} = 1$ (as then $\gamma_1^{2n}$ is the identity). If $\om_1$ has infinite order in $S^1$, then there are infinitely many elements in $\{\om_1^n, n\in \Z-\{0\}\}$, and so there exists some value of $n$ for which $B_{i,j}(\om_1^n, \om_2, a)\neq 0$ for all $i, j\in \{1, 2\}$; this in turn implies that we get our contradiction.
\\\\ Thus, $\om_1$ has order $N$ for some $N\in \Z_{>0}$, and Lemma 3.5 tells us that $N$ is odd. It can be easily checked (by computer, for example) that $B_{1, 2}(\om_1, \om_2, a)= 0$ has no solutions over the $N$th roots of unity, so we are further restricted to only considering $B_{1, 1}, B_{2, 1}$, and $B_{2, 2}$. Now, replacing $\om_1$ with $\om_1^{-1}$ must still leave one of $B_{1, 1}$, $B_{2, 1}$ or $B_{2, 2}$ satisfied (else we are done); this gives us nine pairs of simultaneous equations which all have explicit solutions for $\om_1\in \C$. Checking by computer, we find that the only solutions for $\om_1$ within the odd roots of unity are when $\om_1$ is a cube root of unity. Lemma 3.5 tells us that there exists an extension $c_0$ of $c$ which is $\Gamma$-distinguishing and preservd by $\gamma_1$, so we have now proved Lemma 3.4. $\hfill \blacksquare$
\\\\ Lemma 3.4 gives us very specific information about when precolorings of $\R^2-\{\text{  four points  }\}$ cannot be extended to distinguish $SE(2)$. We can now use this to prove Theorem 3.3.
\\\\ \textbf{Proof of Theorem 3.3.} Let $N>>0$, let $|W|\subset \R^2$ be such that $|W|=N$, and let $c$ be a precoloring of $\R^2-W$. Assume for the sake of contradiction that $c$ cannot be extended to distinguish $SE(2)$. Then, Lemma 3.4 tells us that there exists an extension $c_0$ of $c$ which is $\Gamma$-distinguishing, and invariant under some order 3 rotation $\gamma_0$. Identify a set $Y = \{y_1, y_2, y_3, y_4\}$ of four special points such that $Y$ contains no equilateral triangle. Let $W^\prime = W-Y$. Now, from $W^\prime$ pick a set of four points $S_4 = \{x_4^{(i)}, 1\leq i\leq 4\}$, and three additional points $x_1, x_2, x_3$. There are $F(N)=$ $N-4 \choose 4$$(N-8)(N-9)(N-10)$ such ordered collections $(S_4, x_1, x_2, x_3)$. On the other hand, there are only $O(\f{F(N)}N)$ such ordered collections such that $S_4\cup\{x_1, x_2, x_3\} \cup Y$ contains any equilateral triangle (because for any pair of points $(x_1, x_2)$, there are only two points in $\R^2$ which form an equilateral triangle with them). 
\\\\ For any collection $(S_4, x_1, x_2, x_3)$ such that $S_4\cup\{x_1, x_2, x_3\}\cup Y$ does not contain any equilateral triangle, for each $1\leq j\leq 3$, let $c_j^*$ be the precoloring of $\R^2 -\{y_1, y_2, y_3, y_4\}$ which differs from $c_0$ at exactly $x_j$, and let $c_4^*$ be the precoloring which differs from $c_0$ on exactly $S_4$. By Lemma 3.4, there exists an extension $c_j$ of $c_j^*$ which is $\Gamma$ invariant and symmetric under some rotation $\gamma_j$ of order $3$; we can arrange that the $\gamma_j$ are all of the form $z\mapsto e^{\f{2\pi i}3}z + a_j$. Let $(S_4, x_1^{(1)}, x_2, x_3)$ and $(S_4, x_1^{(2)}, x_2, x_3)$ be two such collections. If $\gamma_1^{(1)} = \gamma_1^{(2)}$, then $\gamma_1^{(1)}$ fixes two colorings $c_1^{(1)}$ and $c_1^{(2)}$ which differ for at least two points ($x_1^{(1)}$ and $x_1^{(2)}$); thus, since $\gamma_1$ has order $3$, $Y\cup \{x_1^{(1)}, x_1^{(2)}\}$ must contain an equilateral triangle. Since $Y\cup \{x_1^{(1)}\}$ does not contain an equilateral triangle, for any fixed $x_1^{(1)}$, there are at most $k$ possible values of $x_1^{(2)}$ for which $\gamma_1^{(1)} = \gamma_1^{(2)}$ could be the case, where $k$ is small and independent of $N$. The same holds true for varying $x_2$ and $x_3$. 
\\\\ We know that $\gamma_4$ has exactly one fixed point; therefore, at least three elements of $S_4$ are not fixed by $\gamma_4$. Furthermore, for each of these three $x_4^{(j)}$, either $\gamma_4(x_4^{(j)})\not\in S_4\cup Y\cup \{x_3\}$ or $\gamma_4^{-1}(x_4^{(j)})\not\in S_4\cup Y\cup \{x_3\}$ (this is because $S_4\cup Y\cup \{x_3\}$ contains no equilateral triangle). Without loss of generality, this means that $\gamma_4^{-1}(x_4^{(j)})\not\in S_4\cup Y\cup \{x_3\}$ for $1\leq j\leq 2$ (two of the three have to satisfy the property for the same power of $\gamma_4$, which can be $\gamma_4^{-1}$ without loss of generality). Then, among these two $x_4^{(j)}$, pick $x_4$ such that $\gamma_4^{-1}(x_4)\neq \gamma_3^{-1}(x_3)$. Thus far, we have simply picked one element of $S_4$, given a fixed collection $(S_4, x_1, x_2, x_3)$.
\\\\ Again because $Y\cup S_4\cup \{x_1, x_2, x_3\}$ contains no equilateral triangle, we may suppose that $\gamma_2^{-1}(x_4)\not\in Y\cup S_3\cup\{x_1, x_2, x_3\}$ (the power of $\gamma_2$ does not matter). We will now find a collection $(S_4, x_1, x_2, x_3)$ which gives us a contradiction from the relation $\gamma_1\gamma_2^{-1}\gamma_3\gamma_4^{-1}=\gamma_3\gamma_4^{-1}\gamma_1\gamma_2^{-1}$ (by construction, $\gamma_1\gamma_2^{-1}$ and $\gamma_3\gamma_4^{-1}$ are translations). 
\\\\ First, note that we have already shown that $\gamma_4^{-1}(x_4)\not\in S_4\cup Y\cup \{x_3\}$, so $c_3(\gamma_3\gamma_4^{-1}(x_4))=c_3(\gamma_4^{-1}(x_4))=c_4(\gamma_4^{-1}(x_4))=c_4(x_4)$. Now $c_3(\gamma_3\gamma_4^{-1}(x_4))=c_2(\gamma_3\gamma_4^{-1}(x_4))$ as long as $\gamma_3\gamma_4^{-1}(x_4)\not\in Y\cup \{x_2, x_3\}$. We already know that $\gamma_3\gamma_4^{-1}(x_4)\neq x_3$. Let $S_4, x_1$, and $x_2$ be fixed. Then, the equations $\gamma_3\gamma_4(x_4) = x^*$ for $x^*\in Y\cup \{x_2\}$ each have exactly one solution for $\gamma_3: z\mapsto \om z + a_3$ (i.e., we can solve the linear equation in $a_3$), so there are at most five such ``bad'' rotations. We also showed that there are at most $k$ choices for $x_3$ yielding any given rotational symmetry, so there are at most $5k$ choices of $x_3$ that are potentially problematic. In total, this means that at most $5k$ $N-4\choose 4$ $(N-8)(N-9) = O(\f{F(N)}N)$ ordered collections which are problematic at this step.
\\\\ Suppose that our collection is not one of those problematic collections. Then, our work so far has shown that $c_2(\gamma_2^{-1}\gamma_3\gamma_4^{-1}(x_4))=c_2(\gamma_3\gamma_4^{-1}(x_4))=c_4(x_4)$. We may change the $c_2$ on the left hand side of this equation to a $c_1$ as long as $\gamma_2^{-1}\gamma_3\gamma_4^{-1}(x_4)\not\in Y\cup\{x_1, x_2\}$. For any $x^*\in Y\cup\{x_1, x_2\}$, the equation $\gamma_3\gamma_4^{-1}(x_4) = \gamma_2(x^*)$ yields at most $6k$ more ``bad'' choices of $x_3$ given a fixed $(S_4, x_1, x_2)$, and again at most $O(\f{F(N)}N)$ problematic collections. 
\\\\ Supposing that our collection is still not problematic, our work so far shows that $c_1(\gamma_1\gamma_2^{-1}\gamma_3\gamma_4^{-1}(x_4))=c_4(x_4)$. To ensure that $\gamma_1\gamma_2^{-1}\gamma_3\gamma_4^{-1}(x_4))\not\in Y\cup\{x_1, x_3\}$, by the same argument an in the previous paragraphs, we have to remove from consideration $O(\f{F(N)}N)$ problematic collections -- the only difference here is that we impose a constraint on the choice of $x_2$ given any fixed $S_4, x_1, x_3$. Thus, for all but $O(\f{F(N)}N)$ of the $F(N)$ possible collections, we find that $c_3(\gamma_1\gamma_2^{-1}\gamma_3\gamma_4^{-1}(x_4))=c_4(x_4)$. 
\\\\ We can use the same exact argument to show that for all but $O(\f{F(N)}N)$ of the $F(N)$ collections, $c_3(\gamma_3\gamma_4^{-1}\gamma_1\gamma_2^{-1}(x_4))=c_2(x_4)$; the only point at which there could be a problem is if $\gamma_2^{-1}(x_4) = x_2$, but we have already arranged for this not to be the case. Therefore, we obtain a contradiction by using any of $F(N)-O(\f{F(N)}N)$ collections $(S_4, x_1, x_2, x_3)$; thus, for sufficiently large $N$, all sets $W$ such that $|W|\geq N$ satisfy $P(W, SE(2))$. This completes the proof of Theorem 3.3. $\hfill \blacksquare$
\section{Infinite extension numbers in higher dimensions}
In [2], it was conjectured that $\ext_D(\R^2, E(2))=7$, $\ext_D(S^2, O(3))= 9$, and $\ext_D(\R^3, E(3)) = 10$. The authors of [2] also posed the question of computing these extension numbers in higher dimensions. In Section $8$, we focused on the $\R^2$ case, where some progress was made; however, the conjecture from [2] remains open. On the other hand, once we go beyond $\R^2$ to $X = \R^n$ ($n\geq 3$) or $X = S^n$ ($n\geq 2$), we now show that the extension number $\ext_D(X, \text{Isom}(X))$ is always infinite (indeed, we will see in Section 5.2 that this is even the case for small subgroups of $\text{Isom}(X)$). However, we separate the case of $S^2$ from the others, because for $X=\R^n$ and $X=S^n$ when $n\geq 3$, we can give very explicit uncountable sets $W$ and precolorings of $X-W$ which cannot be extended.
\subsection{A proof of Theorem 2} 
\textbf{Example 4.1.1}. Let $X = \R^n$ for $n\geq 4$, and let $\{e_1, ... ,  e_n\}$ be the standard basis for $X$. Let $W = \R e_1 \cup \{e_2, e_3, ..., e_n\}$, and note that the pointwise stabilizer of $W$ inside $E(n)$ is trivial, because any isometry which fixes the standard basis of $\R^n$ fixes all of $\R^n$. 
\\\\ \textbf{Claim.} $P(W, E(n))$ does not hold.
\\\\ \textbf{Proof.} Let $c$ be the precoloring of $X-W$ which is uniformly red, and let $c^*$ be any extension of $c$ to $X$. Then, since $|\{e_2, ..., e_n\}|\geq 3$, there exist $e_i$ and $e_j$ with $1\not\in\{i, j\}$ such that $c^*(e_i)=c^*(e_j)$. Furthermore, there exists $\gamma\in O(n)\subset E(n)$ such that $\gamma(e_i)=e_j$, $\gamma(e_j)=e_i$, and $\gamma$ pointwise stabilizes the orthogonal complement of $\R e_i + \R e_j$. Since $e_i$ and $e_j$ are the only possible elements of $\R e_i + \R e_j$ not colored red, we immediately obtain that $\gamma$ fixes $c^*$. $\hfill \blacksquare$
\\\\ \textbf{Example 4.1.2.} Let $X = \R^3$, and let $\{e_1, e_2, e_3\}$ be the standard basis. Considering $\R e_1 + \R e_2 \cong \R^2$, let $T$ be the vertices of an equilateral triangle in $\R e_1 + \R e_2$ centered at the origin. Then, let $W = T \cup \R e_3$. Since $W$ linearly spans all of $\R^3$, it has trivial pointwise stabilizer. 
\\\\ \textbf{Claim.} $P(W, O(3))$ does not hold. 
\\\\ \textbf{Proof.} Let $c$ be the precoloring of $\R^3 - W$ which is uniformly red, and let $c^*$ be any extension of $c$ to $\R^3$. Let $S^1$ be the copy of the unit circle inside $\R^3$ which contains $T$. Since $\text{Isom}(S^1)$ does not distinguish this coloring of $S^1$, there exists some planar reflection or rotation which fixes $c^*|_{S^1}$. Let $\gamma\in O(3)$ act by this reflection or rotation on $\R e_1 + \R e_2$ and stabilize $e_3$ (such an element of $O(3)$ certainly exists). Then, because $\gamma$ stabilizes $\R e_3$ and preserves $c^*|_{S^1}$, since $c^*(\R^3 - \R e_3 \cup S^1) = R$, we conclude that $\gamma$ preserves $c^*$. $\hfill \blacksquare$
\\\\ \textbf{Example 4.1.3.} Let $X = S^n\subset \R^{n+1}$ with $n\geq 3$, and let $e_1, ..., e_{n+1}$ be the standard basis of $\R^{n+1}$. Let $S^1$ be the copy of the unit circle sitting inside $S^n$ with coordinates $x_3 = x_4 = \cdots = x_{n+1} = 0$, and let $T\subset S^1$ be the vertices of an equilateral triangle. Finally, let $S^{n-2}$ be the orthogonal complement of $S^1$ sitting inside $S^n$, and let $W = S^{n-2}\cup T$. Since $W$ linearly spans $\R^{n+1}$, it has trivial pointwise stabilizer within $O(n+1)$. 
\\\\ \textbf{Claim.} $P(W, O(n+1))$ does not hold. 
\\\\ \textbf{Proof.} This is essentially the same as Example 4.1.2. Letting $c$ be the precoloring of $S^n-W$ which is uniformly red, any extension $c^*$ of $c$ to $S^n$, when restricted to $S^1$, is preserved by some planar rotation or reflection. Then, we may let $\gamma\in O(n)$ be the isometry which is equal to this rotation or reflection when restricted to $S^1$ and pointwise stabilizes $S^{n-2}$. By construction, $\gamma$ preserves $c^*$. $\hfill \blacksquare$
\subsection{Extending precolorings on $S^2$: a proof of Theorem 4}
Since all of the counterexamples from Section 5.1 involved invariance under some reflectional symmetry, it is reasonable to ask if removing reflections from the isometry groups would give us finite extension numbers. While one can create counterexamples in $\R^6$ very similar to Example 4.1.1 which do not satisfy $P(W, SO(6))$ (and similarly in higher dimensions), we can also employ another method to create huge numbers of counterexamples in lower dimensions -- even on $S^2$, where Section 5.1 failed to produce any results. 
\\\\ Let us recall the statement of Theorem 4.
\\\\ \textbf{Theorem 4.} \textit{If $W\subset S^2$ is finite, then $P(W, SO(3))$ does not hold. Assuming the axiom of choice, we may replace ``finite'' with ``countable.''}
\\\\ In particular, if we let $SO(3)$ act on $\R^n$ by acting on a particular copy of $\R^3\subset \R^n$, Theorem 4 implies that $\ext_D(\R^n, SO(3))=\infty$. In some sense, this is a much stronger result than Theorem 2, because it produces a huge class of counterexamples. On the other hand, it does not produce any uncountable sets $W$ such that $P(W)$ does not hold. 
\\\\ \textbf{Proof of Theorem 4.} Let $W\subset S^2$ be any finite set, with $|W| = n$. We will construct a precoloring of $S^2-W$ which cannot be extended to distinguish $SO(3)$. First, we'll establish a framework to make the problem easier to think about.
\\\\ Let $\Gamma\subset SO(3)$ be a subgroup with generating set $S$, by which we mean that the elements of $S$ and the elements of $S^{-1}$ together generate $\Gamma$ as a group. Then, we can produce a graph $G(W, \Gamma, S)$ as follows: let $V(G) = \Gamma \cdot W$ and $E(G) = \{(x, y, s)\in V(G)\times V(G)\times S: y = s^{\pm 1} x \text{ for some } s\in S\}$. In other words, each edge is labelled by some element of $S$, and there may be more than one edge connecting two vertices. We say that a $2$-coloring $c$ of $G(W, \Gamma, S)$ is \textit{invariant} under a particular $s\in S$ if all $s$-adjacent vertices (that is, pairs $(x, y)$ such that $(x, y, s)\in E(G)$) have the same color under $c$. 
\\\\ We will construct a ``bad'' precoloring $c$ of $S^2-W$ by picking a group $\Gamma$ (with generating set $S$) such that the graph $G(W, \Gamma, S)$ has a particularly nice structure. In particular, we will use the result, attributed to Hausdorff in [4] (1914), that $SO(3)$ contains a copy of $F_2$, the free group on two letters.
\\\\ There are explicit constructions of free subgroups of $SO(3)$; for example, in [6], it is shown that rotations about the angle $\phi$ in the $x$-$y$ plane and in the $x$-$z$ plane generate a free group provided that $\cos(\phi)\in \Q-\{0, \pm 1, \pm \f 1 2\}$. Furthermore, it is well known that $F_2$ contains as a subgroup $F_{\N}$, the free group on countably infinite many letters (for example, see [5]). Let $\Gamma^\prime = F_\N\subset SO(3)$.
\\\\ \textbf{Claim.} There exists a free subgroup $\Gamma\subset \Gamma^\prime$ with infinite generating set $S$ such that the following two statements are true: (1) the connected components of the elements $w_i\in W$ inside $G(W, \Gamma, S)$ are pairwise disjoint (i.e., there are no paths between elements of $W$), and (2) $G(W, \Gamma, S)$ contains no cycles which contain any element of $W$. 
\\\\ \textbf{Proof.} We first find a subgroup satisfying property (2). Let $S^\prime = \{s_1, s_2, ...\}$ be the free generating set of $\Gamma^\prime$. Let $w_1\in W$. Then, $\stab_{\Gamma^\prime}(w_1)$ is a subgroup of $\Gamma^\prime$ which is also abelian, because $\stab_{SO(3)}(\{w_1\})$ is abelian. We know that abelian subgroups of a free group are isomorphic to $\Z$; this is a special case of the Nielsen-Schreier theorem (which relies on the axiom of choice), but this special case does not rely on the axiom of choice (for example, see [5]). Therefore, $\stab_{\Gamma^\prime}(\{w_1\})$ is generated by a single $\gamma\in \Gamma^\prime$ that can be written as a word in finitely many letters $s_{i_1}, s_{i_2}, ..., s_{i_k}$. Letting $S^\prime_1 = S^\prime - \{s_{i_1}, s_{i_2}, ..., s_{i_k}\}$ and $\Gamma^\prime_1 = \langle S^\prime \rangle$, we see that $\stab_{\Gamma^\prime_1}(\{w_1\})$ is trivial. Repeating this process for each of the elements of $W$, we obtain the subgroup $\Gamma^\prime_n\subset \Gamma^\prime$ with infinite generating set $S^\prime_n$ such that $\stab_{\Gamma^\prime_n}(\{w\})$ is trivial for any $w\in W$ - in other words, $(\Gamma^\prime_n, S^\prime_n)$ satisfies property (2). 
\\\\ To prove the full claim, for any pair $w_i, w_j\in W$, we note that there is at most one element $\gamma\in\Gamma^\prime_n$ such that $\gamma w_i = w_j$; this is because if there were two such elements, $\gamma_1, \gamma_2$, then $\gamma_1\gamma_2^{-1}\in \stab_{\Gamma^\prime_n}(w_j)$, which is trivial. Since there are finitely many elements of $W$, there are finitely many $\gamma$'s in total, each of which is a word in finitely many generators $s_{j_1}, s_{j_2}, ..., s_{j_l}$. Letting $S = S^\prime_n - \{s_{j_1}, ..., s_{j_l}\}$ and $\Gamma = \langle S \rangle$, we see that $(\Gamma, S)$ satisfies both (1) and (2), as desired in the claim. $\hfill \blacksquare$
\\\\ Let $(\Gamma, S)$ be as asserted in the claim. Since $S$ is infinite, let $s_1, ..., s_{2^n}$ be $2^n$ elements of $s$. We will construct a precoloring $c$ of $S^2-W$ such that, enumerating the extensions of $c$ by $c_1, ..., c_{2^n}$, $c_i$ is fixed by $s_i$ for each $1\leq i\leq 2^n$. We will construct $c$ as follows: color $\R^2 - G(W, \Gamma, S)$ red. Thus, it will suffice to show in the end that each of the extensions, when restricted to $G(W, \Gamma, S)$, is invariant under some $s_i$. By property (1), we may consider (and color) the connected components of each $w\in W$ separately. 
\\\\ Enumerate the colorings of $W$ by $c_1^*, c_2^*, ..., c_{2^n}^*$, and let $w\in W$. First, we consider the vertices $x\in V(G(W, \Gamma, S))$ which are adjacent to $w$ -- in other words, $x = s_i^{\pm 1} w$ for some $i\in \N$. If $i>2^n$, then define $c(x) = R$. If $i\leq 2^n$, define $c(x) = c_i^*(w)$. Now, since $G(W, \Gamma, S)$ contains no cycles that contain $w$, for every $y\in \Gamma \{w\}$, there exists a unique ``branch'' $x\in V(G(W, \Gamma, S))$ which is adjacent to $w$ such that every path from $y$ to $w$ has $x$ as its second to last vertex. Then, define $c(y) = c(x)$. Under this construction of $c$, it is clear that $c_i(x) = c_i(w)$ if $x$ is $s_i$-adjacent to $w$. Furthermore, if $y_1$ and $y_2$ are $s_i$-adjacent, then it is clear that $y_1$ and $y_2$ have the same branch $x$, so $c(y_1)=c(x)=c(y_2)$. Thus, the coloring $c_i$ is invariant under $s_i$ for each $1\leq i\leq 2^n$, as desired.
\\\\ For $W$ finite, we managed the above construction without invoking the axiom of choice. For $W$ countably infinite, we can do the same exact construction, but we must use the following fact proved by de Groot and Dekker in [3]: assuming the axiom of choice, $SO(3)$ contains a free group $F$ on uncountably many letters. 
\\\\ We need to use this fact because if $W$ is countably infinite, we may need to remove countably many generators from the generating set of $F$ to remove all of the cycles that contain elements of $W$, and more importantly, there are uncountably many colorings of $W$ which need to satisfy some symmetry. However, replacing ``finite'' with ``countable'' and ``countable'' with ``uncountable'' as necessary, the above arugument will construct a precoloring of $S^2-W$ which cannot be extended to distinguish $F$, or $SO(3)$. $\hfill \blacksquare$
\\\\ Finally, we note that Theorem 4 remains true even if we consider $k$-colorings for $k>2$. In $[2]$, the distinguishing extension number is defined for any $k\geq D_\Gamma(X)$ (rather than just $k = D_\Gamma(X)$), but $\ext_D(S^2, O(3), k) = \infty$ for every $k\geq 2$.  
\section{Open questions}
Of the conjectures and questions posed in [2], the conjecture that $\ext_D(\R^2, E(2))=7$ remains open. We pose a weakened version of this conjecture, as well as another related conjecture.
\\\\ \textbf{Conjecture 5.1.} $\ext_D(\R^2, E(2))<\infty$.
\\\\ \textbf{Conjecture 5.2.} $\ext_D(\R^2, SE(2))=4$.
\\\\ We can also ask if Theorem 2 (which holds for $S^n$ and $\R^n$ for $n\geq 3$) applies to $S^2$.
\\\\ \textbf{Question 5.3.} Do there exist uncountable subsets $W\subset S^2$ such that $P(W, O(3))$ does not hold? Such that $P(W, SO(3))$ does not hold? At least in the case of $SO(3)$, we are inclined to believe that this is not the case. 
\\\\ Finally, the motivation for introducing the distinguishing extension number was to better differentiate group actions on sets -- the distinguishing number $D_\Gamma(X)$ is very often 1, 2, or 3, for example. However, $\ext_D(X, \Gamma)$ cannot differentiate between $O(3)$ acting on $\R^3$ and $O(4)$ acting on $\R^4$, among other things. One possible alternative to the extension number is the following.
\\\\ \textbf{Definition.} The \textit{replacement number} $R(X, \Gamma)$ is the smallest $n\in \N$ such that for every $D_\Gamma(X)$-coloring of $X$, we may replace the colors of at most $n$ points in $X$ to obtain a distinguishing coloring of $X$. 
\\\\ For example, Corollary 17.7 states that $R(S^1, O(2))=3$. It is easy to further establish that $R(\R, E(1))=3$ as well. To make sure that $R(X, \Gamma)$ is not bounded in terms of the distinguishing number (in an obvious way, at least), we note that $R(\R^n, E(n))\geq n$, because any $n-1$ points lie on some hyperplane (so the all-red coloring cannot be fixed using $n-1$ points). Since the ``replacement'' constraint is considerably weaker than the ``extension'' constraint, we are led to the following conjecture.
\\\\ \textbf{Conjecture 5.4.} $R(\R^n, E(n))<\infty$.
\section{Acknowledgements}
This research was conducted under the supervision of Joe Gallian at the University of Minnesota Duluth REU, supported by NSF grant 1358659 and NSA grant H98230-13-1-0273. I would like to thank Joe Gallian, as well as program advisors Noah Arbesfeld, Daniel Kriz, and Adam Hesterberg, for their support throughout the research process. I would also like to thank Aaron Abrams, Xiaoyu He, Brian Lawrence, and David Moulton for their insights and suggestions.

\end{document}